\documentclass[11pt]{article}
\usepackage{amsmath,amsfonts,amsthm,amssymb, graphicx}
\usepackage{subfig}
\usepackage{amsmath}
\usepackage{amsfonts}
\usepackage{amssymb}
\usepackage{graphicx}

\providecommand{\U}[1]{\protect\rule{.1in}{.1in}}

\newtheorem{theorem}{Theorem}

\newtheorem{assumption}{Assumption}

\newtheorem{definition}{Definition}

\newtheorem{lemma}{Lemma}

\newtheorem{proposition}{Proposition}
\newtheorem{remark}{Remark}
\theoremstyle{remark}

\newcommand{\authoronethree}[2][]{\hspace*{9pt}{\small\textrm{\uppercase{#2}},$^{*~***}$ \textit{#1}}\par}
\newcommand{\emailthree}[1]{\footnote{\hspace*{-14pt}$^{***}\,$Email address: #1}\par}

\begin{document}

\title{Perfect Sampling for Infinite Server and Loss Systems}
\author{Jose Blanchet and Jing Dong}
\maketitle

\begin{abstract}
% text of abstract goes here!
We present the first class of perfect sampling (also known as exact simulation) algorithms for the
steady-state distribution of non-Markovian loss networks. We use a variation
of Dominated Coupling From The Past for which we simulate a stationary
infinite server queue backwards in time and analyze the running time in
heavy traffic. In particular, we are able to simulate stationary renewal marked point processes in unbounded regions.
We use the infinite server queue as an upper bound process to simulate loss systems.
The running time analysis of our perfect sampling algorithm for
loss systems is performed in the Quality-Driven (QD) and the
Quality-and-Efficiency-Driven regimes. In both cases, we show that our
algorithm achieves sub-exponential complexity as both the number of servers and
the arrival rate increase. Moreover, in the QD regime, our algorithm
achieves a nearly optimal rate of convergence.
\end{abstract}

%\keywords{perfect sampling, coupling from the past, infinite server queues, loss networks, renewal point processes, many-server asymptotics} % insert keywords separated by a semicolon

%\ams{65C05, 68U20}{60K25} % insert the primary Maths Subject Classification number in the first bracket
         % and the secondary ams number(s) in the second bracket
         % e.g. \ams{60E20}{49G03;49F10}

\section{Introduction} 

\label{intro}

We present the first class of exact simulation algorithms for the steady-state
distribution of non-Markovian loss networks. The running time of our
algorithms is analyzed in the context of many server systems in heavy-traffic;
corresponding both to the so-called Quality-Driven (QD) regime, and the
Quality-and-Efficiency-Driven (QED, also known as Halfin-Whitt) regime. In
both cases, we show that our algorithm achieves sub-exponential complexity as
the number of servers and the arrival rate increase. Moreover, in the QD
regime, our algorithm achieves a nearly optimal rate of convergence. So, more
broadly, our contributions are the first to provide exact simulation
methodology with satisfactory running time analysis in the setting of many
server queues in heavy traffic.

Exact simulation consists in sampling without any bias from the steady-state
distribution of a given ergodic process. Since the inception of Coupling From
The Past (CFTP), the most common exact sampling protocol, proposed in the
ground breaking paper by Propp and Wilson \cite{PropWil:1996}, perfect sampling
(also known as \textquotedblleft exact simulation \textquotedblright%
) has become an important part of stochastic simulation. The majority of the
available exact simulation algorithms for queues involve exponential
distributional assumptions (on service times and/or interarrival times) and
very few of such algorithms are applicable in the context of queueing
networks. None of them, up to date, have been designed and analyzed in the
setting of many server systems in heavy-traffic.

Foss and Tweedie \cite{FosTwe:1998} proved that CFTP can be applied if and
only if the underlying process is uniformly geometric ergodic. Murdoch and
Takahara \cite{MurTak:2006} applied CFTP in the context of queueing models, but
mostly with bounded state space. For instance, they consider loss queues with
renewal arrivals but with bounded service times and in this case, CFTP can be
easily implemented. A variation of CFTP, called Dominated CFTP (DCFTP)\cite{Ken:1998}, 
allows one to apply CFTP-type idea to
obtain unbiased samples from the steady-state distribution of ergodic
processes without requiring uniform ergodicity. A nice summary of DCFTP is
given in \cite{FilHub:2010}. The idea is to
construct a stationary process which suitably dominates the process of
interest and that can be simulated backwards in time from a stationary state
at time zero. Then, a suitable lower bound process, coupled with the upper
bound, must also be simulated in stationarity and backwards in time. A typical
application of DCFTP involves the construction of the upper and lower bound up
to a time in the past when they both meet. Then one says that the coalescence
occurs. The process of interest is reconstructed forward in time from the
coalescence position up to time zero, using the same input sequence that was
used to simulate the upper and lower bounds. The state of the process of
interest at time zero must then follow the corresponding steady-state distribution.

The paper \cite{CorTwe:2001} is one of the earliest to
consider DCFTP in the setting of geometrically ergodic Harris recurrent Markov
chains. General DCFTP algorithms have been developed more recently in \cite{Ken:2004} and \cite{ConKen:2007} for Harris recurrent
chains, although there are important practical limitations as outlined on
p$.788$ in \cite{ConKen:2007}. In particular, their
algorithm assumes that one has analytical access to the transition kernel of
the underlying Markov chain after several transitions. A recent paper by
Sigman \cite{Sig:2011} provides an implementable DCFTP algorithm for
multi-server queues with Poisson arrivals, but the algorithm requires rather
strong conditions on stability; in \cite{Sig:2012} the conditions are
relaxed (also in the setting of Poisson arrivals), using a regenerative
technique but the expected termination time of the algorithm is infinite.

In connection to loss queueing systems, we have already mentioned \cite{MurTak:2006}, 
in the setting of queues with non-Poisson input, but bounded service times. In
\cite{BusGau:2012}, the authors develop a class of CFTP algorithms which
combines aggregation and multiple bounding chains, but exponential service
times and interarrival times are required in their development. The papers \cite{Ken:1998}, \cite{KenMol:2000} and \cite{FerFer:2002} 
are close in spirit to the main ideas of our paper
as we take a point process approach to the problem. However, their approach
requires the use of spatial birth and death processes (generally of poisson
type) as the dominating processes and as pointed out in Section 8 of \cite{BertMol:2002}, the algorithms appear to
significantly increase in complexity as the arrival rate increases.

We provide a practical simulation procedure that works under the assumption of
renewal arrivals (having a finite moment generating function) and service time
distribution with finite mean (although in our running time analysis in heavy
traffic we impose additional moment conditions for service times, but we still
are able to cover distributions such as log-normal, which have been observed
to accurately fit service time distributions in many server applications \cite{BroMan:2002}). 
The performance of our procedures has been
successfully tested numerically in \cite{BlanDon:2012}.

In order to implement our strategy in the setting of loss queues, we simulate
a stationary infinite server queue backwards in time as our dominating
process. A small variation from the standard DCFTP protocol just explained is
that we use the upper bound process itself to detect coalescence, thereby
bypassing the need for a lower bound process and improving the running time of
the algorithm. Basically we detect coalescence over a time interval in which
all customers initially present in the infinite server system leave and no
loss of customers occurs.

We summarize our contributions next:

\begin{itemize}

\item[1)] The design and analysis of the first exact sampling algorithm for the
infinite server queue whose running time is shown to be basically linear in
the arrival rate and thus optimal as the steady state of the infinite server
queue, encoding the remaining service time of each customer, requires on
average a vector which grows linearly in the arrival rate. (See Theorem
\ref{th:inf}.)

\item[2)] The design and analysis of the first exact sampling algorithm for loss
networks under non-Markovian arrivals and a heavy-traffic environment. In the
QD regime, where service utilization is strictly less than 100\%, we show that
our algorithm has near optimal (linear in the arrival rate) running time. In
the QED regime, when the traffic utilization converges to 100\% at a square
root speed as a function of the arrival rate, we show that our algorithm has a
subexponential running time. (See Theorems \ref{th:qd} $\&$ \ref{th:qed}).

\end{itemize}

We point out that our algorithms allow to simulate stationary renewal
processes with independent and identically distributed (i.i.d.) marks in the
positive line on unbounded regions (having finitely many points almost
surely). This connection has been noted in \cite{BlanDon:2012} with a fixed
region (as upposed to a moving frame going backwards in time as we include
here) and without the running time analysis that we perform in this paper for
high arrival rates.

The rest of the paper is organized as follows. Section \ref{sec:main} contains
several subsections. There we introduce our notation and describe the general
strategy to simulate the dominating infinite server process. Then, we describe
how to detect the coalescence time using only the dominating infinite server
process. Once the basic procedures and notation have been explained we proceed
to give precise running time results. The whole description initially
concentrates on a single station and in the end we provide the extension to
the case of loss networks. In Section \ref{sec:alg} we provide the details
required to implement our general strategy outlined in Section \ref{sec:main}
for the infinite server queue in steady state. In Section \ref{sec:per} we
study the running time of our algorithms, some technical results in the
development of this section are given in Appendix.

\section{Basic strategy and main results} \label{sec:main}

In this section we introduce the basic strategy to simulate the systems. We
also present some results about the efficiency of our algorithms. We leave the
details of the algorithms and proofs of the results to subsequent sections. We
start with the strategy of a many-server loss station in steady state and then
generalize our strategy to cover loss networks. 
The method we use is a variation of DCFTP and we use the infinite server queue
as the dominating process.

%We use the infinite server queue
%as the dominating process. The infinite server queue and the loss queue with
%$C$ servers can be coupled naturally: just label the servers in the infinite
%server queue, assign customers to the empty server with the smallest label,
%and by tracking only the state of the first $s$ servers in the infinite server
%queue one automatically tracks the state of the loss queue. So, we can
%easily use the same stream of customers as the input process for both the
%infinite server queue and the many-sever loss queue. And clearly the infinite
%server system serves as an upper bound for the loss system.
%
%As explained in the introduction, we need to simulate the infinite server
%queue backwards in time in stationarity. We then identify a coalescence time,
%that is, a time $T<0$ at which \textit{the state of the loss system is
%completely identified} from the coupled infinite server system. If we know
%such time $T$, we can then simulate the loss queue from $T$ forwards in time
%using the same stream of customers on $[T,0]$ and output the state of the loss
%queue in stationarity at time zero.

To facilitate our explanation, we start with a formal description of the
state of the infinite server (GI/GI/$\infty$) queue. 
%and then identify asuitable coalescence time.

\subsection{Description of the GI/GI/$\infty$ system} \label{sec:not}

%Simulating the infinite server queue in stationarity and backwards in time is
%not trivial, so we first need to explain how to do this task. 

We first introduce some notations and assumptions next. 
Let $N=\{N\left(t\right)  :t\in(-\infty,0]\}$ be a one sided time stationary renewal point
process. We write $\{A_{n}:n\geq1\}$ for the times at which the process $N$
jumps counting backwards in time from time zero with $A_{n+1}<A_{n}<0$.
Furthermore, we define $X_{n}=|A_{n+1}-A_{n}|$. Now let $\{V_{n}:n\geq1\}$ be
a sequence of i.i.d. random variables (r.v.'s) which are independent of the
process $N$. Define $Z_{n}=\left(  A_{n},V_{n}\right)  $ and consider the
marked point process $\mathcal{M}=\{Z_{n}:n\geq1\}\in\mathbb{R}^{2}$ which we
call the \textquotedblleft arriving customer stream\textquotedblright. More
specifically, we consider customers arriving to the system according to a
renewal process with i.i.d. interarrival times $X_{n}$'s. Independent of the
arrival process, their service requirements $V_{n}$'s are also i.i.d..

\begin{figure}[tbh]
\centering
\includegraphics[width=0.5\textwidth]{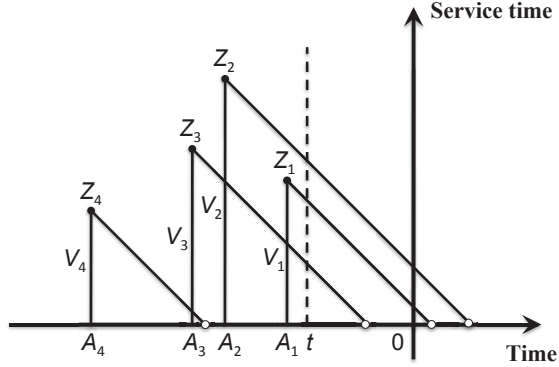} 
\caption{Point process description of an infinite server queue}%
\label{fig:inf}%
\end{figure}

Figure \ref{fig:inf} elaborates on the point process description of the infinite server queue and
is important for describing our simulation strategy.  In Figure \ref{fig:inf}, the point $Z_{n}=(A_{n},V_{n})$ 
denotes the $n$-th customer (counting backward in time), whose arrival time is $A_n$ and service requirement is $V_n$, $n=1,\dots,4$.
One important feature of infinite server queue is that every customer starts service immediately upon arrival (there is no queue). 
If we project $Z_n$ to the horizontal axis by drawing a $-45^{o}$ line. 
The intersection of this line with the horizontal axis is the departure time of such $n$-th
customer.  We follow the technical tradition that an arrival at time $t$ is counted in the system at time $t$ (closed circle) 
while a departure at time $t$ is not counted (open circle).
We can also draw a vertical line at any $t\in\mathbb{R}$. The height of
the intersection of the $-45^{o}$ lines emanating from the points $Z_{n}$ with
$A_{n}\leq t$ and such vertical line, if positive, represents the
corresponding remaining service time of that customer at time $t$.

%Every customer who arrives to the system starts service (i.e. is accepted)
%immediately in the infinite server queue. In pictures (see Figure
%\ref{fig:inf}), the point $Z_{n}=(A_{n},V_{n})$, representing the arrival time
%and service time of the $n$-th customer, is projected to the horizontal axis
%(representing the time axis) by drawing a $-45^{o}$ line. The intersection of
%this line with the horizontal axis is the departure time of such $n$-th
%customer. We can draw a vertical line at any $t\in\mathbb{R}$. The height of
%the intersection of the $-45^{o}$ lines emanating from the points $Z_{n}$ with
%$A_{n}\leq t$ and such vertical line, if positive, represents the
%corresponding remaining service time of that customer at time $t$. \\

%\Fig{Point process description of an infinite server queue} \label{fig:inf}

We write $G(\cdot)=P(X_{n} \leq\cdot)$ for the cumulative distribution
function (CDF) of $X_{n}$ and put $\bar G(\cdot)=1-G(\cdot)$ for its tail CDF.
Similarly, we write $F(\cdot)=P(V_{n} \leq\cdot)$ as the CDF of $V_{n}$ and
$\bar F(\cdot)=1-F(\cdot)$ as its tail CDF.

The following assumption is imposed throughout our discussion:

\begin{assumption} \label{ass:1}
$EX_n<\infty$ and $EV_n<\infty$.
\end{assumption}

We next introduce a Markovian description of the system. Let $Q(t,y)$ denote
the number of people in the system at time $t$ with residual service time
strictly greater than $y$. Notice that for fixed $t$, $Q(t,\cdot)$ is a
piecewise constant step function. If we denote $\{r_{1}(t),..,r_{m}(t)\}$ as
the ordered (positive) remaining service times of customers in the system at
time $t$. Then $Q(t,0)=m$ and $Q(t,y)=\sum_{i=1}^{m}I(r_{i}(t)>y)$. We also
let $E(t)$ denote the time elapsed since the previous arrival at time $t$
(i.e. $E\left(  t\right)  =t-\max\{A_{n}:A_{n}\leq t\}$) and
$W(t)=(E(t),Q(t,\cdot))\in\mathbb{R}^{+}\times\mathcal{D}[0,\infty)$. Then
$\{W(t):t\in R\}$ forms a Markov process which describes the infinite server queue.

Similarly, we denote $W^L(t)=(E^L(t),Q^L(t,\cdot))\in\mathbb{R}^{+}\times\mathcal{D}[0,\infty)$
as the state of the loss queue at time $t$, 
where $E^L(t)$ denote the time elapsed since the previous arrival,  
and $Q^L(t,y)$ counts the number of people in the loss system at time with residual service time
strictly greater than $y$.

Assume both systems start empty from the infinite past, and we use the same stream of customers to update the infinite 
server queue and the loss queue with $C$ servers . 
Then the two systems can be coupled naturally: just label the servers in the infinite
server queue, assign customers to the empty server with the smallest label,
and by tracking only the state of the first $C$ servers in the infinite server
queue one automatically tracks the state of the loss queue.
Based on the above description, we have $E(t)=E^L(t)$ and
$Q^L(t,y) \leq Q(t,y)$ for any $t \in \mathbb{R}$ and $y\geq 0$. Thus
$$W^L(t) \leq W(t)$$
We next define the coalescence time. 
\begin{definition} Coalescence time is a time $T<0$ at which \textit{the state of the loss system is
completely identified} from the coupled infinite server system, i.e. $W^L(t)=W(t)$.
\end{definition}

\subsection{Coalescence time with an $GI/GI/C/C$ queue} \label{sec:coa}

As discussed earlier the infinite server system imposes an upper bound on the
loss system. A natural way to construct the coalescence (or coupling) time
would be to define the coalescence time as the first time (going backwards in
time) the infinite server queue empties (assuming, say, unbounded
interarrival time distribution, this will occur). However, this coalescence
time generally grows exponentially with the arrival rate \cite{Kelly:1991}. 
So, to detect the coalescence in a more efficient
way, we consider the following construction. Let $R(t)$ denote the maximum
remaining service time among all customers in the system at time $t$. And
consider a random time $\tau<0$ satisfying

\begin{itemize}
\item[1)] $R(\tau)<|\tau|$;
\item[2)] $\inf_{\tau \leq t\leq \tau +R(\tau )}\{C-Q(t,0)\}\geq 0$.
\end{itemize}

As we will show in Section \ref{sec:prof:coa}, $\tau$ is well defined and our
coalescence time is $T:=\tau+R(\tau)$. In simple words, since the infinite
sever queue has less than $s$ customers on $[\tau,\tau+R(\tau)]$, the loss
queue is also operating below capacity $C$ on that interval. Everyone who was
present at time $\tau$ in the infinite server queue will have left at time
$\tau+R(\tau)$. Thus the infinite serve queue and the loss queue must have
the same set of customers present in the system by that time. From then on we can recover the state of the loss
queue at time zero using the same stream of customers as for the infinite server queue on
$[\tau+R(\tau),0]$.

\subsection{Basic strategy and main results for the GI/GI/$\infty$ system} \label{sec:inf}

Simulating the infinite server queue in stationarity and backwards in time is
not trivial, so we first need to explain how to do this task. 
There are two cases to be considered. 

\begin{itemize}
\item[Case 1]
The interarrival time has finite exponential moment in a neighborhood of the origin.
More specifically, define $\psi \left( \theta \right) =\log E\exp \left( \theta X_{n}\right) $.
There exists $\theta >0$ such that $\psi (\theta )<\infty $.
\item[Case 2]
The interarrival time does not have finite exponential moment, i.e. it has heavy-tail distribution.
\end{itemize}

%This reduction will allow us
%to concentrate our efforts on the first case. We will come back to the second
%case when we discuss running times.

As we shall explain, we can always reduce the second case to the first one by
defining yet another coupled upper bound process trough truncation.  
Specifically, define $X_{n}\wedge b=\min\{X_{n}
, b\}$. We then fix a suitably large constant $b$ and define a coupled
infinite server queue with truncated interarrival times: $\{X_{n} \wedge b :
n\geq1\}$. This truncation essentially speed up the arrival process. By
coupling we mean we use the same stream of customers to update both the
original system and the truncated one, 
i.e., We use $(X_n, V_n)$ to update the original system and $(X_n\wedge b, V_n)$ to update the truncated one.
We also define the event times as the
arrival time and the departure time of the $n$th customer, $n\geq1$ (counting
backwards in time). Then the infinite server queue with truncated interarrival
times imposes an upper bound, in terms of the number of customers in the
system, on the original infinite server queue at the corresponding event
times, 
i.e. $A_n=\sum_{i=1}^n X_n$ corresponds to $A_n(b):=\sum_{i=1}^n (X_n \wedge b)$ in the truncated system, 
and $A_n+V_n$ corresponds to $A_n^b + V_n$ in the truncated system.
Notice that the actual time of the events (such as arrivals and
departures) may be different for the two systems because of the truncation.
But from the simulation point of view, we simulate the same amount of
information to get the corresponding event times in both systems. In what
follows, we shall first concentrate our discussion on Case 1 which also
includes the infinite server queue with truncated interarrival times. We then
explain how to extend the result to the heavy-tailed case.

\begin{figure}[tbh]
\centering
\includegraphics[width=0.35\textwidth]{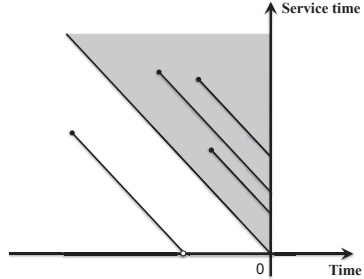} 
\caption{Coupling time of  the infinite server queue }
\label{fig:coup}
\end{figure}

We first introduce the procedure to simulate the state of the stationary
infinite server queue at time zero. We notice from Figure \ref{fig:coup} that
customers $Z_{n}=\{A_{n}, V_{n}\}$, with $V_{n}\leq|A_{n}|$ will have left the
system by time $0$. Thus if we can find a random number $\kappa$ such that
\[
V_{n}\leq|A_{n}| \mbox { for all } n\geq\kappa,
\]
then we can simulate the arrival stream backwards in
time up to $\kappa$ (i.e. $\{Z_{n}:1\leq n\leq\kappa\}$) to recover the state
of the system at time zero. The challenge here is that $\kappa$ defined above
depends on future customer information, i.e. $\{Z_n: n > \kappa\}$. In what follows, we shall explain the elements behind the simulation
of $\kappa$.

We write $\mu=EX_{n}$ and fix an $\epsilon\in(0,\mu)$. Consider any random
number $\kappa$ finite with probability one but large enough such that
\[
A_{n+1}\geq n(\mu-\epsilon)\mbox{ and }V_{n+1}\leq n(\mu-\epsilon) \mbox{ for all } n\geq\kappa.
\]

Let $\kappa(A)$ be a random
time satisfying that $A_{n+1}\geq n(\mu-\epsilon)$ for $n\geq\kappa(A)$, and
$\kappa(V)$ be a random time satisfying that $V_{n+1}\leq n(\mu-\epsilon)$ for
$n\geq\kappa(V)$. Then we can set $\kappa=\max\{\kappa(A),\kappa(V)\}$. The
following proposition states that $\kappa<\infty$ almost surely (a.s.). The
proof is given in Appendix \ref{app:prop}.

%\Fig{Coupling time of infinite server queue } \label{fig:coup}

\begin{proposition} \label{prop}
Under Assumption \ref{ass:1}, the random number $\kappa$ defined above is finite with
probability one.
\end{proposition}

As $\{A_{n}:n\geq1\}$ and $\{V_{n}:n\geq1\}$ are independent of each other, 
the above construction allows us to sample $\{V_{n}:
n\geq1\}$ with $\kappa(V)$, and $\{A_{n}: n\geq1\}$ with $\kappa(A)$ separately.
We next explain the basic sampling strategies.

For the $\{V_n\}$ process, define $J(0):=0$ and $J(l)=\inf\{n > J(l-1): V_{n+1} > n (\mu-\epsilon)\}$ for
$l=1,2,\cdots$. It is the times at which $V_{n+1}$'s exceed the corresponding
increasing boundary $n(\mu-\epsilon)$. Let $\gamma=\inf\{l \geq1:
J(l)=\infty\}$. Then $\kappa(V)=J(\gamma-1)+1$. We first simulate $J(l)$'s for
$l=1,2,...,\gamma-1$, and then simulate the $V_{n}$'s conditional on $J(l)$'s
(see Section \ref{sec:alg:ser} for details).

For the $\{A_n\}$ process, define
\[
\tilde S_{n} = n(\mu-\epsilon)-(A_{n+1}-A_{1})=\sum_{i=1}^{n} Y_{i},
\]
where $Y_{i}=(\mu-\epsilon)-X_{i+1}$. Note that $Y_{i}$'s are i.i.d. with
$EY_{i}=-\epsilon$. Set $\tilde S_{0}=0$. $\{\tilde S_{n}: n \geq0\}$ is a
random walk with negative drift. $A_{n+1}=A_{1}- \tilde S_{n}+n(\mu-\epsilon
)$. If we can simulate some random time $\kappa^{*}$ such that $\tilde S_{n}
\leq0$ for $n \geq\kappa^{*}$, then $|A_{n+1}-A_{1}| \leq n(\mu-\epsilon)$ for
$n \geq\kappa^{*}$. Fix any $m>0$. Define $\Gamma(0)=0$ and $\Delta(l)=\inf\{n
\geq\Gamma(l-1): \tilde S_{n} \leq-m\}$, $\Gamma(l)=\inf\{n \geq\Delta(l):
\tilde S_{n} - \tilde S_{\Delta(l)} \geq m\}$. Let $\alpha=\inf\{l \geq1:
\Gamma(l)=\infty\}$. We notice that $\tilde S_{n} $ will never go above $0$
from $\Delta(\alpha)$ on; which implies that we can set $\kappa(A)=\Delta
(\alpha)$. As we assume the moment generating function of $X_{n}$ is finite in
a neighborhood of the origin, the moment generating function of $Y_{n}$ is
also finite in a neighborhood of the origin. We simulate $\tilde S_{n}$'s
jointly with $\Delta(l)$'s and $\Gamma(l)$'s until $\alpha$ using exponential
tilting and the acceptance rejection method (see Section \ref{sec:alg:arr} for
details).

For the heavy-tailed case (Case 2), we can choose the truncation parameter $b$
such that $E[X_{n}\wedge b]=\int_{0}^{b}\bar G(x)dx=\mu-1/2\epsilon$. This is
doable because we assume $EX_{n}=\int_{0}^{\infty} \bar G(x)dx<\infty$. Set
$\epsilon^{\prime}= 1/2 \epsilon$. Then $E[X_{n} \wedge b]-\epsilon^{\prime
}=\mu-\epsilon$. Denote $A_{n}(b)$ as the backwards renewal times of the
truncated arrival process and let $\kappa(A(b))$ be a random time satisfying
that $|A_{n+1}(b)| \geq n(E[X_{n} \wedge b]-\epsilon^{\prime})$ for $n \geq
\kappa(A(b))$. Then we have $|A_{n+1}| \geq|A_{n+1}(b)| \geq n(\mu-\epsilon)$
for $n \geq\kappa(A(b))$, thus we can set $\kappa(A)=\kappa(A(b))$.

Our algorithm works only under the mild condition in Assumption \ref{ass:1}.
But we do impose stronger conditions on the service time distribution to
rigorously show good algorithmic performance, especially in heavy traffic
(i.e. as the arrival rate increases).

We consider a sequence of systems indexed by $s\in\mathbb{N}^{+}$. We shall
say that $s$\textit{\ is the scale of the system.} 
%The reader might note that
%we have used $s$ in our discussion at the beginning of Section 2 to denote the
%number of servers in the associated loss system. Trying to minimize the
%notational burden we are deliberately avoiding the introduction of new
%notation for the arrival rate parameter. Here we are only considering the
%infinite server system, so in fact, $s$ could take any positive value.
We speed up the arrival rate of the $s$-th system by scale $s$. That is, the
interarrival times of the $s$-th system are given by $X_{n}^{(s)}=X_{n}/s$. We
keep the service time distribution fixed for all systems, i.e. the service
times do not scale with $s$. The following theorem summarizes the performance
of the procedure we proposed for simulating stationary infinite server queue.

\begin{theorem} \label{th:inf}
Assume $E[X_n]<\infty$, and
\begin{itemize}
\item[(1)] if $EV_n^{q}<\infty$ for some $q>2$, then
\begin{equation*}
E_\pi^s\kappa=O(s^{q/(q-1)});
\end{equation*}
\item[(2)] if we further assume $E[\exp(\theta V_n)]<\infty$ for some $%
\theta>0$, then
\begin{equation*}
E_\pi^s \kappa=O(s \log s).
\end{equation*}
\end{itemize}
\end{theorem}

We prove it by establishing two bounds for $\kappa(A)$ and $\kappa(V)$ respectively. The details is given in Section \ref{sec:prof:inf}.

We next extend the procedure to simulate states of the stationary infinite
server system backwards in time for time intervals of any
specified length. The construction is very similar to the single time point (i.e.
time zero) case explained above.

Define $\kappa_{0}:=1$. We consider a sequence of random times $\kappa_{j}$,
$j=1,2, \cdots$, finite with probability one but large enough such that
\begin{equation}
\label{con:term}|A_{n}-A_{\kappa_{j-1}}|\geq(n-\kappa_{j-1})(\mu
-\epsilon)\text{ and }V_{n}\leq(n-\kappa_{j-1})(\mu-\epsilon)
\mbox{ for all } n\geq\kappa_{j}.
\end{equation}
Notice that $V_{n}\leq|A_{n}-A_{\kappa
_{j-1}}|$ for $n\geq\kappa_{j}$. This implies that a customer who arrives
before $A_{\kappa_{j}}$ will not be in the system at time $A_{\kappa_{j-1}}$.
Thus, using $\{Z_{n}: 1\leq n\leq\kappa_{j}\}$, we can
recover the system descriptor $W(t)$ for $t\in\lbrack A_{\kappa_{j-1}},0]$.

%Figure \ref{fig:cons} gives more details about the construction. Every point
%$Z_{n}$, with $n>\kappa_{j}$, will not land into the upper triangle defined by
%the vertical line at $A_{\kappa_{j-1}}$ and the $-45^{o}$ line intersecting it
%at the time axis (x axis).
%
%%\Fig{Coupling time of infinite server queue} \label{fig:cons}
%
%\begin{figure}[tbh]
%\centering
%\includegraphics[width=0.5\textwidth]{fig3.eps}
%\caption{Coupling times of the infinite server queue}%
%\label{fig:cons}%
%\end{figure}

The $\kappa_{j}$'s give us some flexibility to separate the simulation of the
two processes. We first simulate the service times and then conditional on the
sample path of the service time we simulate the arrival process jointly with
$\kappa_{j}$'s.

Define $J_{1}(0):=1$ and let
\begin{align*}
J_{k}(l)  &  = \inf\{ n > J_{k}(l-1): V_{n} >(n-J_{k}(0))(\mu-\epsilon)),\\
\gamma_{k}  &  = \inf\{l \geq0: J_{k}(l)=\infty\},\\
J_{k+1}(0)  &  = J_{k}(\gamma_{k}-1)
\end{align*}
for $k=1,2,\cdots$ and $l=1,2,\cdots,\gamma_{k}$.\newline 
We first simulate the
random time: $J_{k}(l)$'s for $k=1,2,\cdots$ and $l=1,2,\cdots,\gamma_{k}$, and
then simulate $\{V_{n}: n \geq1\}$ conditional on $J_{k}(l)$'s; see Algorithm
I in Section \ref{sec:alg:ser} for details.

Given the sample path of $\{V_{n}:n\geq1\}$ and $J_{k}(l)$'s, we next simulate
$\{A_{n}:n\geq1\}$ and $\kappa_{j}$'s. This is done by simulating the
negative-drift random walk $\tilde{S}_{n}$ jointly with its running time maximum. Define
$\Delta_{1}(0):=0$ and $\Gamma_{1}(0):=0$. Fix $m>0$ and let
\begin{align*}
\Delta_{j}(l)  &  =\inf\{n\geq\Gamma_{j}(l-1):\tilde{S}_{n}-\tilde{S}%
_{\Delta_{j}(0)}\leq-m\},\\
\Gamma_{j}(l)  &  =\inf\{n\geq\Delta_{j}(l):\tilde{S}_{n}-\tilde{S}%
_{\Delta_{j}(l)}\geq m\},\\
\alpha_{j}  &  =\inf\{l\geq1:\Gamma_{j}(l)=\infty\},\\
\kappa_{j}  &  =\min\{J_{k}(0):J_{k}(0)\geq\Delta_{j}(\alpha_{j})+1\},\\
\Delta_{j+1}(0)  &  =\kappa_{j}-1,\\
\Gamma_{j+1}(0)  &  =\Delta_{j+1}(0)
\end{align*}
for $j=1,2,...$ and $l=1,2,\cdots,\alpha_{j}$.\newline Notice that the process
$\tilde{S}_{n}$ will never go above $\tilde{S}_{\Delta_{j}(0)}$ from
$\Delta_{j}(\alpha_{j})$ on. This implies that $|A_{n}-A_{\kappa_{j-1}}%
|\geq(n-\kappa_{j-1})(\mu-\epsilon)$ for $n\geq\kappa_{j}$. 
Under the light-tail assumption (Case 1), we simulate the random
times $\Delta_{j}(l)$ and $\Gamma_{j}(l)$ for $j=1,2,...$, $l=1,2,...,\alpha
_{j}$ and $\{\tilde{S}_{n}:n\geq0\}$ by the exponential tilting and
acceptance-rejection method. The details are explained in Algorithm II in
Section \ref{sec:alg:arr}.

For the heavy-tailed case (Case 2), we again simulate the infinite server
queue with truncated interarrival times first. We carefully choose the
truncation parameter $b$ and $\epsilon^{\prime}$ such $E[X_{n}\wedge
b]-\epsilon^{\prime}$ coincides with $\mu-\epsilon$. Then the $\kappa_{j}%
(b)$'s we constructed for the truncated system must automatically satisfy the
conditions characterizing $\kappa_{j}$'s in (\ref{con:term}) for the original
system as well.

\subsection{Basic strategy and main results for the $GI/GI/C/C$ system}

\label{sec:loss}

Once we simulate the customer streams backwards in time and construct the
states of the dominating stationary infinite server queue accordingly, we can
check and find the coalescence time $T=\tau+R(\tau)$ where $\tau$ is defined
in Section \ref{sec:coa} backwards in time. Use the state of the infinite
server queue at time $T$ as the state of the many-server loss queue at the
same time and go forwards in time using the same stream of customers to
construct the state of the loss queue up to time $0$.

Like in the infinite server queue case, we again consider a sequence of
systems indexed by $s\in\mathbb{N}^+$ where the arrival rate of the $s$-th
system is scaled by $s$ and the service rate is kept fixed. Let $\rho
=E[V_{n}]/E[X_{n}]$ (the ratio of the mean service time and mean interarrival
time of the base system). We analyze the system in two heavy-traffic
asymptotic regimes. One is the quality driven (QD) regime where $\rho<1$ and
the number of servers in the $s$-th system, $C_s$, is $s$. The other is the quality
and efficiency driven (QED) regime where $\rho=1$ and the number of servers in
the $s$-th system, $C_s$, is $s+b\sqrt{s}$ with $b>0$. 

Theorem \ref{th:qd} summarizes the performance of the coalescence time in the
QD regime.

\begin{theorem}\label{th:qd}
Assume $EX_n<\infty$ and $X_n$'s are non-lattice and strictly positive. We also assume that $EV_n^q<\infty$ for any $q>0$ and the
cumulative distribution function (CDF) of $V_n$ is continuous. Then
\begin{equation*}
E_\pi^s \tau = o(s^{\delta})
\end{equation*}
for any $\delta>0$.
\end{theorem}

\begin{remark}
The existence of all moments assumption on the service time distribution
covers a range of heavy tailed distributions, such as Weibull and
log-normal, which are known to fit well data in applications \cite{BroMan:2002}.
\end{remark}

Theorem \ref{th:qed} analyze the performance of the coalescence time in the
QED regime. 

\begin{theorem} \label{th:qed}
Assume $EX_n^2<\infty$.
We also assume $EV_n^q<\infty$ for any $q>0$ and the CDF of $V_n$ is continuous.
Then for $b$ large enough, we have
\begin{equation*}
\log E_{\pi}^s\tau=o(s^\delta)
\end{equation*}
for any $\delta>0$.
\end{theorem}

The main difficulty in the proof of Theorem \ref{th:qd} and Theorem \ref{th:qed} is that it involves the state of the system on an interval rather than a single point.
In Section \ref{sec:prof:coa}, we prove Theorem \ref{th:qd} by using the sample path large deviation results \cite{BlanChenLam:2012} of infinite server queue.  For Thereom \ref{th:qed}, we prove it by applying Borel-TIS inequality \cite{Adl:1990} to  the diffusion limit process of infinite server queue \cite{WhitPang:2010}. The details is also given in Section \ref{sec:prof:coa}. 

\subsection{Extensions and main results for the loss network}

\label{sec:net}

Following the definition in \cite{Kelly:1991}, we consider a generalized loss network with $J$ stations, labeled
$1,2,\cdots,J$ and suppose that station $j$ comprises $C_{j}$ servers. We have
$L$ possible routes, labeled $1,2,...,L$ and for each route $l$, a $J$
dimensional routing vector $P_{l}$. $P_{l}$ is consist of $1$'s and $0$'s,
where $P_{l}(j)=1$ means route $l$ requires a server at station $j$. A routing
request $l$ is blocked and thus lost if any station $j$ with $P_{l}(j)=1$ is
full at the arrival time of the request. Customers requesting route $l$ form a
renewal process with i.i.d. interarrival times $\{X_{n}^{(l)}:n\geq1\}$. The
CDF of $X_{n}^{(l)}$ is $G_{l}$. Independent of the arrival process, the
service times $\{V_{n}^{(l)}:n\geq1\}$ are also i.i.d. with CDF $F_{l}$. We
assume that $G_{l}$'s and $F_{l}$'s satisfy Assumption \ref{ass:1}.

Following the same strategy as in the many-server loss queue case, we first
couple the loss network with a network of infinite-server stations. Notice
that no customer is blocked or lost in the infinite server system, thus it
imposes an upper bound on the number of jobs in the loss system. Let
$Q_{j}(t,y)$ denote the number of jobs in the $j$-th station with remaining
service time strictly greater than $y$ at time $t$. Note that a class $l$ job with
remaining service time greater than $y$ in the system will be counted in
all $Q_{j}(t,y)$'s with $P_l(j)=1$. Let
$R_{j}(t)$ denote the longest remaining service time among all customers in
station $j$ at time $t$. Let $R(t)=\max_{1\leq j\leq J}\{R_{j}(t)\}$. Then
similar to the many server loss queue, we define a random time $\tau^{\prime}$
satisfying the following conditions:

\begin{itemize}
\item[1)] $R(\tau^{\prime }) \leq |\tau^{\prime }|$,
\item[2)] $\inf_{\tau^{\prime }\leq t \leq \tau^{\prime}+R(\tau^{\prime })}\inf_{1 \leq
j \leq J}\{C_j-Q_j(t,0)\} \geq 0$,\newline
i.e. all links are operating below capacity on the interval $[\tau^{\prime
}, \tau^{\prime}+R^{m}(\tau^{\prime })]$.
\end{itemize}

At time $\tau^{\prime}+R(\tau^{\prime})$, everyone in the network of
infinite-server stations will be in the loss network as well. Thus from then
on (forwards in time), we can update the loss system using the inputs of the
infinite-server system.

In order to simulate the network of infinite-server stations with $L$ types of
routing requests, we simulate $L$ independent networks of infinite-server
stations; each dealing with a single type of routing request. Then we do a
superposition of them. The simulation of each independent network of
infinite-server stations are exactly the same as what we have described in
Section 2.3, as a type $l$ routing request occupies a server from each station
$j$ with $P_{l}(j)=1$ simultaneously and for the same amount of time. For the
$l$-th system, let $Z_{n}^{(l)}=(A_{n}^{(l)},V_{n}^{(l)})$ represent the
arrival time and service time of the $n$-th routing request counting backwards
in time and $\kappa^{(l)}$ be a random time satisfying that $V_{n}^{(l)}%
\leq|A_{n}^{(l)}|$ for all $n\geq\kappa^{(l)}$. Then following the procedure
described in Section 2.3, we will be able to simulate $\kappa^{(l)}$ as the
maximum of two random times associated the arrival process and service time
process respectively.

We now consider a sequence of systems indexed by $s\in\mathbb{N}^{+}$. We
speed up the the arrival rate of the $s$-th system by $s$, i.e. $X_{n}%
^{(l,s)}=X_{n}^{(l)}/s$, and keep the service rate fixed. The same result as
in Theorem \ref{th:inf} will still be holding here. Specifically,

\begin{theorem}[Theorem \ref{th:inf}'] \label{th:inf:1}
Assume $EX_n^{(l)}<\infty$ (C).
\begin{itemize}
\item[(1)] if $E[(V_n^{(l)})^{q}]<\infty$ for some $q>2$, then
\begin{equation*}
E_\pi^s\kappa^{(l)}=O(s^{q/(q-1)});
\end{equation*}
\item[(2)] if we further assume $E[\exp(\theta V_n^{(l)})]<\infty$ for some $%
\theta>0$, then
\begin{equation*}
E_\pi^s \kappa^{(l)}=O(s \log s)
\end{equation*}
for $l=1,2,\cdots,L.$
\end{itemize}
\end{theorem}

The proof of Theorem
\ref{th:inf:1} is the same as that of Theorem \ref{th:inf} except for a few
notational changes, thus we shall omit it here.

If we held the number of routing request types, $L$, fixed, as we shall
explain below, similar results as in Theorem \ref{th:qd} and Theorem
\ref{th:qed} for the coalescence time will be holding here as well. We
again run $L$ independent networks of infinite-server stations as described
above. Network $l$ serves routing request of type $l$ only, for $l=1,2,...,L$.
Let $Q^{(l)}(t,0)$ denote the number of jobs in network $l$ at time $t$ and
$R^{(l)}(t)$ denote the maximum remaining service time among all jobs in the
network at time $t$. Then we have $R(t)=\max\{R^{(l)}(t):1\leq l\leq L\}$.

We consider two asymptotic regimes. One is the QD regime where for the base
system we have
\begin{equation}
\label{eq:nqd}\sum_{l=1}^{L}\frac{EV_{n}^{(l)}}{EX_{n}^{(l)}}P_{j}(l)<C_{j}.
\end{equation}
For the $s$-th system, the number of servers in the $j$-th station is $C_{j}%
^{s}=sC_{j}$ for $j=1,2,...,J$. 

Assign a fixed number $H_{l}$ to each
route $l$. $H_{l}$ is well chosen such that $E[V_{n}^{(l)}]/E[X_{n}%
^{(l)}]<H_{l}$ and $\sum_{l=1}^{L}H_{l}P_{l}(j)\leq C_{j}$. This is doable
because of (\ref{eq:nqd}). Let $H_{l}^{s}=sH_{l}$. Define a random time
$\bar{\tau}^{\prime}$ satisfying the following two conditions:

\begin{itemize}
\item[1)] $R^{(l)}(\bar \tau^{\prime }) \leq |\bar \tau^{\prime }|$ for $%
l=1,2,\cdots,L$,
\item[2)] $\inf_{\bar\tau^{\prime}\leq t \leq \bar\tau^{\prime}+ R(\bar\tau^{\prime })}
\{H_l-Q^{l}(t,0)\} \geq 0$ for $l=1,2,\cdots,L$.
\end{itemize}

Notice that $\bar\tau^{\prime}$ is an upper bound on $\tau^{\prime}$. As the
number of types of routing request is fixed at $L$ (it does not scale with
$s$), using the construction outlined in Section
\ref{sec:prof:qd}, we can show that the result in Theorem \ref{th:qd} holds for
$\bar\tau^{\prime}$ as well.

\begin{theorem}[Theorem \ref{th:qd}']\label{th:qd:1}
Assume $EX_n^{(l)}<\infty$ and $X_n^{(l)}$'s are non-lattice and strictly positive.
We also assume $E[(V_n^{(l)})^{q}] <\infty$ for any $q>0$ and $F_l$ is
continuous. Then
\begin{equation*}
E_{\pi}^s \tau^{\prime}=o(s^\delta)
\end{equation*}
for any $\delta>0$.
\end{theorem}

The other asymptotic regime is the QED regime where for the base system we
have
\[
\sum_{l=1}^{L} \frac{EV_{n}^{(l)}}{EX_{n}^{(l)}} P_{j}(l) = C_{j}%
\]
and the number of servers in the $j$-th station of the $s$-th system is
$C_{j}^{s}=sC_{j} + \beta_{j} \sqrt s $ for $j=1,2,\cdots,J$

We then let $I_{l} = E[V_{n}^{(l)}]/E[X_{n}^{(l)}]$ and $I_{l}^{s}%
=sI_{l}+a_{l}\sqrt s$ where $a_{l}$'s are well chosen such that $\sum
_{l=1}^{L}a_{l}P_{j}(l)\leq \beta_{j}$.

We define a random time $\tilde\tau^{\prime}$ that satisfies the following two
conditions:

\begin{itemize}
\item[1)] $R^{(l)}(\tilde \tau^{\prime}) \leq |\tilde \tau^{\prime}|$ for $%
l=1,2,\cdots,L$,
\item[2)] $\inf_{\tilde\tau^{\prime}\leq t \leq
\tilde\tau^{\prime}+R(\tilde\tau^{\prime})} \{I_l-Q^{(l)}(t,0)\} \geq 0$ for
$l=1,2,\cdots,L$.
\end{itemize}

As before, $\tilde\tau^{\prime}$ is an upper bound on $\tau^{\prime}$. It is
easy to check using the construction outlined in Section \ref{sec:prof:qed}
that the result in Theorem \ref{th:qed} holds for $\tilde\tau^{\prime}$ as well.

\begin{theorem}[Theorem \ref{th:qed}']\label{th:qed:1}
Assume $E[(X_n^{(l)})^2]<\infty$.
We also assume $E[(V_{n}^{(l)})^{q}]<\infty $ for any $q >0$.
Then for $b_{j}$'s large enough, we have
\begin{equation*}
\log E_{\pi }^{s}\tau ^{\prime }=o(s^{\delta })
\end{equation*}
for any $\delta >0$.
\end{theorem}

We shall omit the proof of Theorem \ref{th:qd:1} and Theorem \ref{th:qed:1} as
it is the same as the proof of Theorem \ref{th:qd} and Theorem \ref{th:qed}
with the introduction $\bar\tau^{\prime}$ and $\tilde\tau^{\prime}$ except for
a few notational changes.

\section{Detailed simulation algorithms}

\label{sec:alg}

In order to provide the details of our simulation algorithms outlined in
Section 2.3, we shall first work under the light-tailed case (Case 1) where we
assume there exists $\theta>0$ such that $\psi(\theta)<\infty$. The extension
to the heavy-tailed case (Case 2) was introduced in Section
\ref{sec:main} and we shall provide more details in Section
\ref{sec:alg:heavy}.

We further impose the following assumptions on our ability to simulate the
service times and interarrival times.

\begin{assumption} \label{ass:2}
We assume that $F(\cdot)$ is known and easily accessible either
in closed form or via efficient numerical procedures. Moreover, we can
simulate $V_n$ conditional on $V_n \in (a,b]$ with $P(V_n \in (a,b])>0$.
\end{assumption}

\begin{assumption} \label{ass:3}
Suppose that $G\left( \cdot \right) $ is known and that it is
possible to simulate from $G_{eq}\left( \cdot \right) :=\mu ^{-1}\int_{\cdot
}^{\infty }\overline{G}\left( t\right) dt$. Moreover, let $G_{\theta }\left(
\cdot \right) =E\exp (\theta X_{n}-\psi \left( \theta \right) )I(X_{n}\leq
\cdot )$ be the associated exponentially tilted distribution with parameter $%
\theta $ for $\psi \left( \theta \right) <\infty $. We assume that we can
simulate from $G_{\theta }\left( \cdot \right) $.
\end{assumption}

We next introduce our algorithm to simulate $\{V_{n}: n \geq1\}$.
Conditional on the sample path of $\{V_{n}: n \geq1\}$, we then explain how to to simulate $\{A_{n}: n \geq1\}$
and $\kappa_{j}$'s.

\subsection{Simulation of $\{V_{n}: n \geq1\}$ and $J_{k}(l)$'s for
$k=1,2,\cdots$, $l=1,2,\cdots,\gamma_{k}$}

\label{sec:alg:ser}

We will first introduce the procedure to simulate $J_{1}(l)$ for
$l=1,2,\cdots,\gamma_{1}$. Recall that $J_{1}(0):=1$. Let $p(n)=P(V_{1}>n(\mu-\epsilon))$. 
Then $P(J_1(l)=\infty|J_1(l-1)=k)=\prod_{n=k+1}^{\infty}(1-p(n))$. It involves the evaluation
of the product of infinite terms. In Procedure A, we introduce a sandwiching approximation scheme
to accomplish that. 

The following lemma guarantees the termination of our procedure.

\begin{lemma}
If $EV_1 <\infty$, then
\begin{equation}  \label{eq: low}
P(J_1(1)=\infty)=\prod_{n=1}^{\infty}(1-p(n)) \geq \exp(-\frac{cEV_1}{%
\mu-\epsilon})>0
\end{equation}
for some constant $c$ depending on the value of $p(1)$, and consequently 
$$E\gamma_1\leq \exp(cEV_1/(\mu-\epsilon))<\infty.$$
\end{lemma}

\begin{proof}
\begin{align*}
P(J_{1}(1)=\infty) = \prod_{n=1}^{\infty}(1-p(n))  &  \geq\prod_{n=1}^{\infty}
\exp(-cp(n))\\
&  \geq\exp(-\frac{c}{\mu-\epsilon}\int_{0}^{\infty} P(V_{1}>\nu)d\nu) =
\exp(-\frac{cEV_{1}}{\mu-\epsilon}).
\end{align*}
For $l=2,3,\cdots$, conditional on $J_{1}(l-1)=k$:
\begin{align*}
P(J_{1}(l)=\infty|J_{1}(l-1)=k)  &  = \prod_{n=k+1}^{\infty}(1-p(n))\\
&  \geq\exp(-\frac{c\int_{k}^{\infty} P(V_{1}>\nu)d\nu}{\mu-\epsilon})\geq
\exp(-\frac{cEV_{1}}{\mu-\epsilon}),
\end{align*}
thus $\gamma_{1}$ is stochastically dominated by a geometric random variable
with parameter $p=\exp(-cEV_{1}/(\mu-\epsilon))$. The result then follows.
\end{proof}

We next introduce our sandwiching approximation scheme.
Notice that
\begin{equation}
\label{eq: two}\prod_{i=k+1}^{h}(1-p(i))\geq P(J_{1}(l)=\infty| J_{1}%
(l-1)=k)\geq\prod_{i=k+1}^{h}(1-p(i))\times\exp(-\frac{2\int_{h}^{\infty}
P(V_{1}>\nu)d\nu}{\mu-\epsilon})
\end{equation}
for $h > k$.\newline Another important observation is that if we let
$\prod_{i=k+1}^{k} (1-p(i))=1$,
\[
\prod_{i=k+1}^{h-1}(1-p(i))-\prod_{i=k}^{h}(1-p(i))=p(h)\prod_{i=k}%
^{h-1}(1-p(i))=P(J_{1}(l)=h|J_{1}(l-1)=k)
\]
for $h > k.$\newline Let
\[
u(h)= \exp(-\frac{2\int_{h}^{\infty} P(V_{1}>\nu)d\nu}{\mu-\epsilon}).
\]
We now propose the following procedure to simulate the value of $J_{1}(l)$
conditional on $J_{1}(l-1)=k.$ \\

\noindent\textbf{Procedure A (Simulate }$J_{1}(l)$\textbf{ given }%
$J_{1}(l-1)=k$\textbf{)}

\begin{enumerate}
\item Initialize $h=k+1$, $g=1-p(h)$ and $f=gu(h)$. Simulate $U \sim$ Unif$[0,1]$
\item While $f < U < g$,\\
set $h=h+1$, $g=g(1-p(h))$ and $f=gu(h)$\\
end while
\item If $U\leq f$, then $J_1(l)=\infty$. Otherwise, $J_1(l)=h$.\\
\end{enumerate}

The simulation of $J_{k}(l)$ for $l=1,2,...,\gamma_{k}$ follows the same
rationale. We let $p_{k}(n)=P(V_{1}>n(\mu-\epsilon)| V_{1} \leq(n+J_{k}%
(0)-J_{k-1}(0))(\mu-\epsilon))$. Then following the same argument leading to
(\ref{eq: low}) and (\ref{eq: two}), we have correspondingly
\[
P(J_{k}(1) = \infty) > 0,
\]
and
\begin{align*}
\prod_{i=n+1}^{h}(1-p_{k}(i)) 
&  \geq P(J_{k}(l)-J_{k}(0)=\infty|
J_{k}(l-1)-J_{k}(0)=n)\\
&  \geq\prod_{i=n+1}^{h}(1-p_{k}(i))\times\exp(-\frac{2\int_{h}^{\infty}
P(V_{1}>\nu| V_{1} \leq\nu+(J_{k}(0)-J_{k-1}(0))(\mu-\epsilon))d\nu}%
{\mu-\epsilon})
\end{align*}
for $h > n$.\newline Let
\[
u_{k}(h)=\exp(-\frac{2\int_{h}^{\infty} P(V_{1}>\nu| V_{1} \leq\nu
+(J_{k}(0)-J_{k-1}(0))(\mu-\epsilon))d\nu}{\mu-\epsilon}).
\]
We now propose a modification of Procedure A that allows us to simulate
$J_{k}(l)$ conditional on $J_{k}(l-1)-J_{k}(0)=n$.\newline

\noindent\textbf{Procedure A1 (Simulate }$J_{k}(l)$\textbf{ given }%
$J_{k}(l-1)-J_{k}(0)=n$\textbf{)}

\begin{enumerate}
\item Initialize $h=n+1$, $g=1-p_k(h)$ and $f=gu_k(h)$. Simulate $U \sim$
Unif$[0,1]$.
\item While $f < U < g$,\\
set $h=h+1$, $g=g(1-p_k(h))$ and $f=gu_k(h)$\\
end while
\item If $U\leq f$, then $J_k(l)=\infty$. Otherwise, $J_k(l)=J_k(0)+h$.
\end{enumerate}

Based on Procedure A1 and our previous analysis we have:\newline

\noindent\textbf{Algorithm I} (Sample $V_{n}$'s jointly with $J_{k}(l)$'s)

\begin{itemize}
\item[Step 0.] Set $J_0(0)=-\infty$, $J_1(0)=1$, $k=1$, $l=1$. Simulate $V_1$
according to its nominal distribution.
\item[Step 1.] Simulate $J_k(l)$ conditional on the value of $J_k(l-1)$
using Procedure A1.
\item[Step 2.] If $J_k(l)=\infty$, set $\gamma_k=l$, $J_{k+1}(0)=J_k(%
\gamma_k-1)$, $k=k+1$, $l=1$ and go back to Step 1. Otherwise, go to Step 3.
\item[Step 3.] Simulate $V_n$ for $J_k(l-1)< n < J_k(l)$ by conditioning on $%
V_n \leq (n-J_k(0))(\mu-\epsilon)$ and simulate $V_{J_k(l)}$ by conditioning
on $(J_k(l)-J_k(0))(\mu-\epsilon)< V_{J_k(l)} \leq
(J_k(l)-J_{k-1}(0))(\mu-\epsilon)$. Set $l=l+1$ and go back to Step 1.
\end{itemize}

When running the above algorithm, we specify $K$ as the number of intervals
($[J_{k}(0), J_{k}(\gamma_{k}-1)]$) we want to simulate. We then run Algorithm
I from $k=1$ till $k=K$. The program will give us $\{V_{n}: 1 \leq n \leq
J_{K}(\gamma_{K}-1)\}$ and $J_{k}(l)$'s for $k=1,2,\cdots,K$, $l=1,2,\cdots
,\gamma_{k}$.

\subsection{Simulation of $\{A_{n}: n \geq1\}$ and $\Delta_{j}(l)$'s,
$\Gamma_{j}(l)$'s for $j=1,2,...$, $l=1,2,...,\alpha_{j}$}

\label{sec:alg:arr}

Given the sample path of $\{V_{n}: n \geq1\}$, we will first explain how to
simulate the $\Delta_{j}(l)$'s and $\Gamma_{j}(l)$'s sequentially and jointly
with the underlying random walk $\{\tilde S_{n}: n \geq1\}$. We then simulate
$A_{1}$ according to $G_{eq}(\cdot)$ and set $A_{n+1}=A_{1}+n(\epsilon
-\mu)-\tilde S_{n}$. The analysis and methodology in this subsection closely
follows those in \cite{EnsGlyn:2000} and \cite{BlanSig:2011}. 
The same procedure can be used to simulate a negative drifted random walk, $\tilde S_n$,
together with its running time maximum defined as $\max_{k\geq n}\{\tilde S_k\}$.

Let $\mathcal{F}_{n}=\sigma\{Y_{1},Y_{2},\cdots,Y_{n}\}$, the $\sigma$-field
generated by the $Y_{j}$'s up to time $n$. Let $\xi\geq0$ and define
\[
T_{\xi}:=\inf\{n\geq0:\tilde{S}_{n}>\xi\}.
\]
Then by the strong Markov property we have that for $1\leq l\leq\alpha_{j}$,
\[
P(\Gamma_{j}(l)=\infty|\mathcal{F}_{\Delta_{j}(l)})=P(\Gamma_{j}%
(l)=\infty|\tilde{S}_{\Delta_{j}(l)})=P(T_{m}=\infty)>0,
\]
where we use $P\left(  \cdot\right)  $ to denote the nominal probability measure.

It is important then to notice that
\[
P(\alpha_{j}=k) =P(T_{m}<\infty)^{k-1}P(T_{m}=\infty)
\]
for $k\geq1$. In other words, $\alpha_{j}$ is geometrically distributed. The
procedure that we have in mind is to simulate each stage $\Delta_{j}%
(\alpha_{j})$ in time intervals, and the number of time intervals is precisely
$\alpha_{j}$.

Let $\psi_{Y}(\theta)=\log E\exp(\theta Y_{i})$ be the log moment generating
function of $Y_{i}$. As we assume $\psi_{X}(\theta)$ is finite in a
neighborhood of zero, $\psi_{Y}(\cdot)$ is also finite in a neighborhood of
zero. Moreover $EY_{i}=\psi_{Y}^{\prime}(0)=-\epsilon$ and Var$(Y_{i}%
)=\psi_{Y}^{\prime\prime}(0)>0$. Then by the convexity of $\psi_{Y}(\cdot)$,
one can always select $\epsilon>0$ sufficiently small so that there exists
$\eta>0$ with $\psi_{Y}(\eta)=0$ and $\psi_{Y}^{\prime}(\eta) \in(0, \infty)$.
The root $\eta$ allows us to define a new measure $P_{\eta}$ based on
exponential tilting so that
\[
\frac{dP_{\eta}}{dP}\left(  Y_{i}\right)  =\exp(\eta Y_{i}).
\]
Moreover, under $P_{\eta}$, $\tilde{S}_{n}$ is random walk with positive drift
equal to $\psi_{Y}^{\prime}\left(  \eta\right)  $ \cite{Asm:2003}. Therefore $P_{\eta}(T_{\xi}<\infty)=1$ and
\[
q(\xi) :=P(T_{\xi}<\infty)=E_{\eta}\exp(-\eta\tilde{S}_{T_{\xi}})
\]
for each $\xi\geq0$. Based on the above analysis we now introduce a convenient
representation to simulate a Bernoulli random variable $J\left(  \xi\right)  $
with parameter $q\left(  \xi\right)  $, namely,
\begin{equation}
J\left(  \xi\right)  =I(U\leq\exp(-\eta\tilde{S}_{T_{\xi}})), \label{ConRep}%
\end{equation}
where $U$ is a uniform random variable independent of everything else under
$P_{\eta}$.

Identity (\ref{ConRep}) provides the basis for an implementable algorithm to
simulate a Bernoulli random variable with success probability $q(\xi)$.
Sampling $\{\tilde{S}_{1},\cdots,\tilde{S}_{T_{\xi}}\}$ conditional on
$T_{\xi}<\infty$, as we shall explain now, corresponds to basically the same
procedure. First, let us write
\[
P^{\ast}(\cdot)=P(\cdot|T_{\xi}<\infty).
\]
The following result provides an expression for the likelihood ratio between
$P^{\ast}$ and $P_{\eta}$.

\begin{lemma}
We have that
\begin{equation*}
\frac{dP^{\ast}}{dP_{\eta}}(\tilde{S}_{1},...,\tilde{S}_{T_{\xi}})=\frac {%
\exp(-\eta\tilde{S}_{T_{\xi}})}{P(T_{\xi}<\infty)}\leq\frac{\exp(-\eta \xi)}{%
P(T_{\xi}<\infty)}.
\end{equation*}
\end{lemma}

\begin{proof}
\begin{align*}
P(\tilde S_{1} \in H_{1}, ... , \tilde S_{T_{\xi}} \in H_{T_{\xi}} | T_{\xi} <
\infty)  &  = \frac{P(\tilde S_{1} \in H_{1}, ... , \tilde S_{T_{\xi}} \in
H_{T_{\xi}}, T_{\xi} < \infty)}{P(T_{\xi} < \infty)}\\
&  =\frac{E_{\eta}[\exp(-\eta\tilde S_{T_{\xi}})I(\tilde S_{1} \in H_{0}, ...
, \tilde S_{T_{\xi}} \in H_{T_{\xi}})]}{P(T_{\xi}< \infty)}.
\end{align*}
\end{proof}

The previous lemma provides the basis for a simple acceptance / rejection
procedure to simulate $\{\tilde{S}_{1},...,\tilde{S}_{T_{\xi}}\}$ conditional
on $T_{\xi}<\infty$. More precisely, we propose $\{\tilde{S}_{1},...,\tilde
{S}_{T_{\xi}}\}$ from $P_{\eta}\left(  \cdot\right)  $. Then one generates a
uniform random variable $U$ independent of everything else and accept the
proposal if
\[
U\leq\frac{P(T_{\xi}<\infty)}{\exp(-\eta\xi)}\times\frac{dP^{\ast}}{dP_{\eta}%
}(\tilde{S}_{1},...,\tilde{S}_{T_{\xi}})=\exp(-\eta(\tilde{S}_{T_{\xi}}%
-\xi)).
\]
This criterion coincides with $J(\xi)$ according to (\ref{ConRep}). So, the
procedure above simultaneously obtains both a Bernoulli r.v. $J(\xi)$ with
parameter $q(\xi)$, and the corresponding path $\{\tilde{S}_{1},...,\tilde
{S}_{T_{\xi}}\}$ conditional on $T_{\xi}<\infty$ under $P(\cdot)$ if
$J(\xi)=1$.

As $E[Y_{i}]=-\epsilon<0$, by strong law of large numbers we have $\Delta
_{j}(l)<\infty$ almost surely for $j=1,2,...$ and $l=1,2,...,\alpha_{j}$. We
next define
\[
\bar q(\xi): =1-q(\xi)=P(T_{\xi}=\infty)
\]
and
\[
P^{\prime}(\cdot)=P(\cdot|T_{\xi}=\infty).
\]
The following result provides an expression for the likelihood ratio between
$P^{\prime}$ and $P$.

\begin{lemma}
We have that
\begin{equation*}
\frac{dP^{\prime}}{dP}(\tilde{S}_{1},...,\tilde{S}_{n})=\frac{%
I(T_{\xi}>l)\bar q(\xi-\tilde{S}_{n})}{P(T_{\xi}=\infty)}\leq\frac{1}{%
P(T_{\xi}=\infty)}.
\end{equation*}
\end{lemma}

\begin{proof}
\begin{align*}
& P(\tilde S_{1} \in H_{1}, .... , \tilde S_{n} \in H_{n} | T_{\xi}= \infty)\\  
& = \frac{P(\tilde S_{1} \in H_{1}, ... \tilde S_{n} \in H_{n}, T_{\xi}=
\infty)}{P(T_{\xi}= \infty)}\\
&  = \frac{E[I(\tilde S_{1} \in H_{1}, ... , \tilde S_{n}\in H_{n}) I(T_{\xi}>
n) P(T_{\xi}=\infty| \tilde S_{1},...,\tilde S_{n}) ]}{P(T_{\xi}=\infty)}.
\end{align*}
The result then follows from the strong Markov property and homogeneity of the
random walk. 
\end{proof}

We are in good shape now to apply acceptance / rejection to sample from
$P^{\prime}$. The previous lemma indicates that to sample $\{\tilde{S}%
_{1},...,\tilde{S}_{n}\}$ given $T_{\xi}=\infty$. We can propose from the
original (nominal) distribution and accept with probability $\bar q(\xi
-\tilde{S}_{n})$ as long as $\tilde{S}_{j}\leq\xi$ for all $0\leq j\leq n$.
And in order to perform the acceptance test we need to sample a Bernoulli with
parameter $\bar q(\xi-\tilde{S}_{n})$, but this is easily done using identity
(\ref{ConRep}).

Now consider $0 \leq\xi_{1}<\xi_{2}$, we define
\[
P^{o}(\cdot|T_{\xi_{1}}<\infty, T_{\xi_{2}}=\infty).
\]
The following result provides an expression for the likelihood ratio between
$P^{o}$ and $P_{\eta}$.

\begin{lemma}
We have that
\begin{equation*}
\frac{dP^o}{dP_\eta} (\tilde S_1, ... \tilde S_{T_{\xi_1}})=\frac{\exp(-\eta
\tilde S_{T_{\xi_1}})\bar q(\xi_2- \tilde S_{T_{\xi_1}})}{%
P(T_{\xi_1}<\infty, T_{\xi_2}=\infty)}\leq\frac{\exp(-\eta \xi_1)}{%
P(T_{\xi_1}<\infty, T_{\xi_2}=\infty)}.
\end{equation*}
\end{lemma}

\begin{proof}
\begin{align*}
&  P(\tilde S_{1} \in H_{1}, ... , \tilde S_{T_{\xi_{1}}} \in H_{T_{\xi_{1}}}
|T_{\xi_{1}} < \infty, T_{\xi_{2}}=\infty)\\
&  = \frac{E_{\eta}[I(\tilde S_{1} \in H_{1}, ... , \tilde S_{T_{\xi_{1}}} \in
H_{T_{\xi_{1}}})\exp(-\eta\tilde S_{T_{\xi_{1}}})P(T_{\xi_{2}}=\infty| \tilde
S_{1}, ... ,\tilde S_{T_{\xi_{1}}})]}{P(T_{\xi_{1}}<\infty, T_{\xi_{2}}%
=\infty)}.
\end{align*}
\end{proof}

We again use acceptance/rejection to sample $\{\tilde S_{1}, ... , \tilde
S_{T_{\xi_{1}}}\}$ given $T_{\xi_{1}}<\infty$ and $T_{\xi_{2}}=\infty$. We
propose $\{\tilde S_{1}, ..., \tilde S_{T_{\xi_{1}}}\}$ from $P_{\eta}(\cdot
)$. Then we simulate a uniform random variable $U$ independent of all else and
accept the proposal if
\[
U \leq\frac{P(T_{\xi_{1}}<\infty, T_{\xi_{2}}=\infty)}{\exp(-\eta\xi_{1}%
)}\times\frac{dP^{o}}{dP_{\eta}}(\tilde S_{1}, ..., \tilde S_{T_{\xi_{1}}%
})=\exp(-\eta(\tilde S_{T_{\xi_{1}}}-\xi_{1}))q(\xi_{2}-\tilde S_{T_{\xi_{1}}%
}).
\]

Based on the above analysis we propose the following algorithm.\\

\noindent\textbf{Algorithm II} (Given $V_{n}$'s and $J_{k}(l)$'s, sample
$\tilde S_{n}$'s together with $\Delta_{j}(l)$'s, $\Gamma_{j}(l)$'s and
$\kappa_{j}$'s)

\begin{itemize}
\item[Step 0.] Set $\Delta_1(0)=\Gamma_1(0)=0$, $\tilde S_0=0$, $j=1$, $l=1$%
, $\xi=\infty$, $\gamma=-m$. Sample $A_1$ according to $G_{eq}(\cdot)$.
\item[Step 1.] Simulate $S_1, ... , S_{T_\gamma}$ from the original
(nominal) distribution.
\item[Step 2.] If $S_i \leq \xi$ for all $1 \leq i \leq T_\gamma$ then
sample a Bernoulli $J(\xi-S_{T_\gamma})$ with parameter $q(\xi-S_{T_\gamma})$
using (\ref{ConRep}) and continue to step 3. Otherwise (i.e. $S_i > \xi$ for
some $1\leq i\leq T_\gamma$) go back to step 1.
\item[Step 3.] If $J(\xi-S_{T_\gamma})=1$, go back to step 1. Otherwise $%
J(\xi-S_{T_\gamma})=0$, let $\Delta_j(l)=\Gamma_j(l-1)+T_\gamma$ and $\tilde
S_{\Gamma_j(l-1)+i}= \tilde S_{\Gamma_j(l-1)}+S_i$ for $i=1,...,T_\gamma$.
If $j \geq 2$, set $\xi=\tilde S_{\Delta_{j-1}(\alpha_{j-1})}+m-\tilde
S_{\Delta_j(l)}$.
\item[Step 4.] Simulate $S_1,..., S_{T_m}$ from $P_\eta(\cdot)$. Sample a
Bernoulli $J(\xi-S_{T_m})$ with parameter $q(\xi-S_{T_m})$ using (\ref%
{ConRep}) and $U\sim$ Unif$[0,1]$. Let $J^*=I(U \leq
\exp(-\eta(S_{T_m}-m))\times(1-J(\xi- S_{T_m}))$.
\item[Step 5.] If $J^*=1$, let $\Gamma_j(l)=\Delta_j(l)+T_m$ and $\tilde
S_{\Delta_j(l)+i}=\tilde S_{\Delta_j(l)}+S_i$ for $1\leq i\leq T_m$. Set $%
\gamma=\min\{0,\tilde S_{\Delta_j(0)}-m-\tilde S_{\Gamma_j(l)}\}$. If $j
\geq 2$, set $\xi=\tilde S_{\Delta_{j-1}(\alpha_{j-1})}+m-S_{\Gamma_j(l)}$.
Set $l=l+1$ and go back to step 1. Otherwise $J^*=0$, set $\alpha_j=l$, $%
\kappa_j=\inf \{J_k(0): J_k(0) \geq \Delta_{j}(\alpha_{j})+1\}$, $%
\Delta_{j+1}(0)=\kappa_j-1$, $\xi=m$ and continue to step 6.
\item[Step 6.] Let $h=\Delta_{j+1}(0)-\Delta_{j}(\alpha_{j})$. Sample $S_1,
... , S_h$ from the original distribution.
\item[Step 7.] If $S_i \leq \xi$ for all $1 \leq i \leq h$ then sample a
Bernoulli $J(\xi-S_h)$ with parameter $q(\xi-S_h)$ using (\ref{ConRep}) and
continue to step 8. Otherwise (i.e. $S_i>\xi$ for some $1 \leq i \leq h$),
go back to step 6.
\item[Step 8.] If $J(\xi-S_h)=1$, go back to step 6. Otherwise $J(\xi-S_h)=0$%
, let $\tilde S_{\Delta_{j}(\alpha_{j})+i}=\tilde
S_{\Delta_{j}(\alpha_{j})}+S_i$ for $i=1,..., h$. Set $A_{n+1}=A_1+n(%
\epsilon-\mu)-\tilde S_n$ for $n=\Delta_{j}(0)+1,...,\Delta_{j+1}(0)$. Set $%
j=j+1$, $l=1$, $\xi=\tilde S_{\Delta_{j-1}(\alpha_{j-1})}+m-\tilde
S_{\Delta_j(0)}$, $\gamma=-m$ and go back to step 1.
\end{itemize}

When running the above algorithm, we specify $K$ as the number of intervals
($[\kappa_{j-1}, \kappa_{j}]$) we want to simulate and then repeat the above
process from $j=1$ till $j=K$. The program will give us $\{A_{n} : 1 \leq n
\leq\kappa_{K}\}$ and $\{\kappa_{j}: 1 \leq j \leq K\}$.

\subsection{Coupled infinite server queue with truncated interarrival times}

\label{sec:alg:heavy}

In this subsection, we provide some additional details for simulating the
coupled truncated system together with the original system.

We first explain how to simulate $A_{1}$ jointly with $A_{1}(b)$. The
equilibrium distribution of $X_{n}$ is $G_{eq}(x)=\int_{0}^{x} \bar G(u)du /
EX_{n}$ and the equilibrium distribution of $X_{n}\wedge b$ is
\[
G_{eq}^{b}(x)=\frac{\int_{0}^{x}\bar G(u)du}{E[X_{n}\wedge b]} I\{x \leq b\}.
\]
Thus we simulate $A_{1}$ with CDF $G_{eq}(x)$, if $A_{1} \leq b$, we set
$A_{1}(b)=A_{1}$. Otherwise if $A_{1}>b$, we keep simulating $X_{e}$ with CDF
$G_{eq}(x)$ until $X_{e}\leq b$ and set $A_{1}(b)=X_{e}$. In particular we
have $A_{1}(b)\leq A_{1}$.

When simulating $X_{n}\wedge b$'s from the nominal distribution, we first
simulate $X_{n}$ with CDF $G(\cdot)$ and set $X_{n}\wedge b = \min\{X_{n},
b\}$. Denote $Y_{n}(b)=(E[X_{n}\wedge b]- \epsilon^{\prime})- X_{n} \wedge b$
and let $\eta_{b}$ be chosen such that $\log E \exp(\eta_{b} Y_{n}(b))=0$.
When simulating $X_{n}\wedge b$'s under exponential tilting $P_{\eta_{b}%
}(\cdot)$, we first simulate $Y_{n}(b)$ under $P_{\eta_{b}}(\cdot)$ and set
$X_{n} \wedge b=(E[X_{n}\wedge b]-\epsilon^{\prime})-Y_{n}(b)$. If $X_{n}
\wedge b < b$, set $X_{n} =X_{n} \wedge b$, otherwise ($X_{n}\wedge b=b$),
sample $X_{n}$ conditional on $X_{n} \geq b$ under the nominal distribution
$P(\cdot)$.

\section{Performance analysis}

\label{sec:per}

In the previous section, we provide our simulation algorithm and show that our
algorithm works in the sense that the termination time is finite with
probability one. In this section, we conduct some further asymptotic analysis
on the performance of our algorithm. We first analyze the algorithm for the
infinite server system and then conduct some analysis on the coalescence time
for the many-server loss system.

\subsection{Termination time for the infinite server system (Proof of Theorem
\ref{th:inf})}

\label{sec:prof:inf}

Theorem \ref{th:inf} provides the relationship between the moment of the
service times and $E_{\pi}^{s} \kappa$. We next give a proof of it. We shall
omit the subscription $\pi$ and $s$ when there is no confusion for notational
convinience. We first give a proof of the light tailed case. 
Recall that
$\kappa=\max\{\kappa(V), \kappa(A)\}$,
where
$\kappa(V)=\inf\{k > 1: V_{n+1} \leq n(\mu-\epsilon)/s \mbox{ for all } n \geq
k\}$
and
$\kappa(A)= \inf\{k>1: A_{n+1} \geq n(\mu-\epsilon)/s \mbox{ for all } n \geq
k\}$.
We prove the theorem by establishing the bounds for $\kappa(V)$ (Lemma \ref{lm:ser}) and $\kappa(A)$ (Lemma
\ref{lm:arr}) respectively.

\begin{lemma} \label{lm:ser}
If $EV_n^{q}<\infty$ for some $q>2$, then
\begin{equation*}
E\kappa(V)=O(s^{q/(q-1)}).
\end{equation*}
\end{lemma}

\begin{proof} Let $p(n)=P(V_{1} > n (\mu-\epsilon)/s)$. 
For $k$ sufficiently large, we have
\begin{align*}
P(\kappa(V)>k)  &  = 1- \prod_{n=k+1}^{\infty} (1-p(n))\\
&  \leq1-\exp(-\frac{2s}{\mu-\epsilon}\int_{k(\mu-\epsilon)/s}^{\infty}
P(V>\nu)d\nu).
\end{align*}
By Chebyshev's inequality
\[
P(V_{n}>\nu)\leq\frac{EV_{n}^{q}}{\nu^{q}}.
\]
Let $\delta=1/(q-1)$, then for $s$ sufficiently large, we have
\begin{align*}
\sum_{k=s^{1+\delta}}^{\infty} P(\kappa(V)>k)  &  \leq\sum_{k=s^{1+\delta}%
}^{\infty} \frac{2s}{\mu-\epsilon}\int_{k(\mu-\epsilon)/s}^{\infty}
P(V>\nu)d\nu\\
&  \leq\frac{2EV_{n}^{q}s^{q}}{(q-1)(q-2)(\mu
-\epsilon)^{q}} \sum_{k=s^{1+\delta}}^{\infty} \frac{1}{k^{q-1}}\\
&  = O(s^{q-(1+\delta)(\delta-2)}).
\end{align*}
As $q-(1+\delta)(q-2) = 1+\delta$,
\begin{align*}
E\kappa(V)  &  = \sum_{k=0}^{\infty} P(\kappa(V)>k)\\
&  = \sum_{k=0}^{s^{1+\delta}-1} P(\kappa(V)>k) + \sum_{k=s^{1+\delta}%
}^{\infty} P(\kappa(V)>k)\\
&  \leq s^{1+\delta}+O(s^{1+\delta}).
\end{align*}
\end{proof}

Notice that when $E\exp(\theta V_{n})<\infty$ for some $\theta>0$,
\[
P(V_{n}>\nu)\leq E\exp(\theta(V_{n}-\nu)) = E\exp(\theta V_{n})\exp(-\theta
\nu).
\]
Similarly as above, for $s$ sufficiently large we have
\[
\sum_{k= \lceil\frac{2}{ \theta(\mu-\epsilon)} s \log s\rceil}^{\infty}
P(\kappa(V)>k) \leq\frac{2E\exp(\theta V_{n})}{(\mu-\epsilon)^{2}\theta^{2}}%
\]
and
\[
E\kappa(V) = \sum_{k=0}^{s\log s-1} P(\kappa(V)>k) + \sum_{k=s\log s}^{\infty}
P(\kappa(V)>k) \leq s\log s + O(1).
\]
Thus if $E\exp(\theta V)<\infty$ for some $\theta>0$, then
\[
E\kappa(V)=O(s\log s).
\]

\begin{lemma} \label{lm:arr}
Assume there exist $\theta>0$, such that $\psi(\theta) <\infty$, then
\begin{equation*}
E\kappa(A)=O(s).
\end{equation*}
\end{lemma}

\begin{proof} Based on the algorithm proposed in Section \ref{sec:alg:arr},
we divide the proof into two parts. We first prove that the expected number of
iterations is $O(1)$. We then prove that the expected number of steps to reach
$-m$ or $m$ is $O(s)$.

Let $T_{\xi}= \inf\{n \geq0: \tilde S_{n} > \xi\}$. Recall that for the base
system there exist $\eta>0$ with $\psi_{Y}(\eta)=0$ and $\psi_{Y}^{\prime
}(\eta)>0$. And the number of iterations is distributed as a geometric random
variable with probability of success $P(T_{m}=\infty)=1-E_{\eta}\exp
(-\eta\tilde S_{T_{m}})$

Then for the $s$th system with $Y_{i}^{s}=Y_{i}/s$ we have $\tilde S_{n}/s>m$
is equivalent to $\tilde S_{n}>sm$. Thus the number of iterations is a
Geometric random variable with probability of success
\[
P(T_{sm}=\infty)=1- E_{\eta}\exp(-\eta\tilde S_{T_{sm}}) \geq1-\exp(-\eta
sm).
\]

Similarly, let $T_{\xi}^{\prime}= \inf\{n \geq0: \tilde S_{n} < \xi\}$. Define
$M_{n} = \tilde S_{n} + n\epsilon$, then $M_{n}$ is a martingale with respect
to the filtration generated by $\{Y_{1}, Y_{2}, ... , Y_{n}\}$. As
$EY_{i}=-\epsilon<0$, $P(T_{-m}<\infty)=1$. By the Optional Sampling Theorem,
$E M_{T_{-m}^{\prime}}=E \tilde S_{T_{-m}^{\prime}} + \epsilon E
T_{-m}^{\prime}=0$. Thus
\[
ET_{-m}^{\prime}=\frac{m}{\epsilon}-\frac{E[m-S_{T_{-m}^{\prime}}]}{\epsilon
}.
\]
Then for the $s$th system we have
\[
ET_{-sm}^{\prime}=\frac{sm}{\epsilon}-\frac{E[sm-S_{T_{-sm}^{\prime}}%
]}{\epsilon}.
\]
$(sm-S_{T^{\prime}_{-sm}})$ converges to the ladder hight $Y^{-}$ distribution
as $s \rightarrow\infty$ and $\sup_{m} E[(sm-S_{T_{-sm}^{\prime}})^{p}%
]<\infty$ yields $E[(Y^{-})^{p}]<\infty$ for $p>1$ \cite{Asm:2003}. Therefore,
\[
ET_{-sm}^{\prime}=O(s).
\]
\end{proof}

For the heavy-tailed case, we select the truncation parameter $b$ such that
$E[X_{n}\wedge b]=\mu-1/2\epsilon$. Then we set $\epsilon^{\prime}%
=1/2\epsilon$ and define $\kappa(A(b))$ as a random time satisfying that
$|A_{n+1}|\geq n(E[X_{n} \wedge b]-\epsilon^{\prime})=n(\mu-\epsilon)$ for
$n\geq\kappa(A(b))$. As $|A_{n+1}|\geq|A_{n+1}(b)|$ under our coupling scheme,
we can set $\kappa(A)=\kappa(A(b))$. By Lemma \ref{lm:arr}, we have
$E\kappa(A)=E\kappa(A(b))=O(s)$. $\kappa(V)$ is defined as before, a random
time satisfying that $V_{n} \leq n(\mu-\epsilon)$ for $n \geq\kappa(V)$. Then
$E\kappa(V)=O(s\log s)$ by Lemma \ref{lm:ser}.

As $\kappa=(\kappa(V),\kappa(A(b))$, we have $E\kappa=O(s\log s)$. This
concludes the proof of Theorem \ref{th:inf}.

\subsection{Coalescence time for the many-server loss system (Proof of Theorem
\ref{th:qd} and Theorem \ref{th:qed})}

\label{sec:prof:coa}

As we are simulating the process backwards in time, it is natural to define
the following filtration
\[
\overleftarrow{\mathcal{H}}_{t}=\sigma\{W(-u): 0 \leq u \leq t\},
\]
for which $\overleftarrow{\mathcal{H}}_{u} \subset\overleftarrow{\mathcal{H}%
}_{t}$ for $0 \leq u \leq t$. $\tau$ is a stopping time with respect to
$\overleftarrow{\mathcal{H}}_{t}$. We next try to draw connections between the
backward process and some forward process. Define
\[
\tau^{*}=\inf\{t+R(t): \sup_{t \leq u \leq t+R(t)}\{Q(u,0)\} < s, t \geq0\}
\]
$\tau^{*}$ is a stopping time with respect to $\mathcal{H}_{t}$ where
$\mathcal{H}_{t}=\sigma\{M(u): 0 \leq u \leq t\}$. The stochastic process
$\{Q(t,0): t \in{\mathbb{R}}\}$ has a piecewise constant sample path with a
finite number of points of discontinuity on any finite length intervals almost
surely. Thus for any fixed $T>0$, we have
\begin{align*}
P_{\pi}(\tau> T)  &  = P_{\pi}(\bigcap_{-T \leq t \leq0}(\{R(t) > -t\}
\bigcup( \bigcup_{t \leq u \leq(t+R(t))\wedge0} (\{Q(u,0)>s\})))\\
&  = P_{\pi}(\bigcap_{-T \leq t \leq0}(\{R(T+t)>-t\} \bigcup(\bigcup_{T+t \leq
u \leq(T+t+R(T+t))\wedge T}\{Q(u,0)>s\})))\\
&  = P_{\pi}(\bigcap_{0 \leq w \leq T} ( \{R(w)>T-w\} \bigcup(\bigcup_{w \leq
u \leq(w+R(w))\wedge T}\{Q(u,0)>s\})))\\
&  = P_{\pi}(\tau^{*} > T)
\end{align*}
The second equality holds by stationarity; this gives us $E_{\pi}\tau=E_{\pi
}\tau^{*}$. Next, we use a special construction similar to that in Section 4
of \cite{BlanLam:2012} to prove the results for $E_{\pi}^{s}\tau^{*}$.  
The idea is to use a geometric trial argument. We divide the time frame into blocks
that are roughly independent. And if the process is well-behaved 
(staying around its measure-valued fluid limit) on one block, then $\tau^*$ is reached before the end of that block.

Let $\bar Q(t,y)$ denote the number of customers in the infinite server queue
that starts empty at time zero with remaining service time greater than $y$ at
time $t \geq0$. For convenience, we also define $\bar Q_{u}(t,y)=\bar Q(u+t,y)-\bar Q(u,t+y)$ for
$u\leq t$, as the number of customers who arrive after $u$ with remaining
service time larger than $y$ at time $u+t$.

\subsubsection{Proof of Theorem \ref{th:qd}.}

\label{sec:prof:qd}
We first prove the theorem for the light-tailed case. The heavy-tail case
proceeds by selecting the truncation parameter $c$ sufficiently large.

For the QD regime, by ``well-behaved", we mean that the process does not deviate $\delta s$, for some $\delta>0$, from its fluid limit.
The following lemma states that the probability of not being well-behaved decays exponentially fast with the system scale.
 
\begin{lemma}\label{lm:dp1}
Assume $\psi(\theta)<\infty$ for some $\theta>0$ and $X_n$'s are non-lattice and strictly positive. We also assume the CDF
of $V_n$ is continuous. Then for any $\delta>0$, there exist $I^*(\delta)>0$%
, such that
\begin{equation*}
P(\bar Q(t,y) > (1+\delta)\lambda s\int_{y}^{t+y}\bar F(u)du
\mbox{ for some } t\in[0,1], y \in [0,\infty)) = \exp(-sI^*(\delta)+o(s)).
\end{equation*}
\end{lemma}

The proof of Lemma \ref{lm:dp1} follows form the tow-parameter sample path
large deviation result for infinite server queues in \cite{BlanChenLam:2012}. 
We shall omit it here.

We next introduce our construction of ``blocks".
Let $l(s)=\inf\{y: (1+\delta)s\int_{y}^{\infty}\bar F(u)du \leq\frac{1}{2}\}$,
we define the following sequence of random times $\Xi_{i}$'s: $\Xi_{0}:=0$.
Given $\Xi_{i-1}$ for $i=1,2,\cdots$, define
\begin{align*}
r_{i}  &  =\inf\{k: k \geq R(\Xi_{i-1}), k=1,2,\cdots),\\
z  &  = \inf\{k: k \geq l(s), k=1,2,\cdots\},\\
\Xi_{i}  &  = \Xi_{i-1}+r_{i}+z.
\end{align*}
We define a Bernoulli random
variable $\xi_{i}$, with $\xi_{i}=1$ if and only if
\[
\bar Q_{\Xi_{i-1}+(k-1)t_{0}}(t,y) \leq(1+\delta)\lambda s \int_{y}^{t+y}\bar
F(u) du
\]
for all $t \in[0,1]$, $y \in[0,\infty)$ and every $k=1,2,\cdots, r_{i}+z$.  

Choose $\delta< 1/\rho-1 $. We first check
that $\xi_{i}=1$ implies that $\tau^*$ is reached before $\Xi_{i}$.
Since $r_{i} \geq R(\Xi_{i-1})$, all the customers in the system at time
$\Xi_{i-1}+r_{i}$ will be those who arrive after $\Xi_{i}$. Then $\xi_{i}=1$
implies that
\begin{eqnarray*}
Q(\Xi_{i-1}+r_{i}, y) &\leq& \sum_{k=1}^{r_{i}/t_{0}} \int_{(k-1)t_{0}+y}^{kt_{0}+y}\bar F(u)du\\ 
&=& (1+\delta)\lambda s\int_{y}^{r_{i}+y} \bar F(u)du\\
&\leq& (1+\delta)\lambda s \int_{y}^{\infty}\bar F(u)du,
\end{eqnarray*}
thus $R(\Xi_{i-1}+r_{i}) \leq l(s)$.\newline And for every $t\in(k-1,k]$,
$k=1,2,...,z$
\begin{align*}
Q(\Xi_{i-1}+r_{i}+t,y)  &  \leq(1+\delta) \lambda s \int_{y}^{r_{i}+t+y} \bar
F(u) du\\
&  \leq(1+\delta) \lambda s \int_{y}^{\infty} \bar F(u) du,
\end{align*}
thus $Q(\Xi_{i-1}+r_{i}+t,0)\leq(1+\delta)\rho s \leq s$ for $t \in[0,
R(\Xi_{i-1}+r_{i})]$.\newline Now let $N=\inf\{i\geq1: \xi_{i}=1\}$, then
$$E \tau^{*} \leq E\sum_{i=1}^{N}(r_{i}+z).$$

We now show a bound for $E\sum_{i=1}^{N}(r_{i}+z)$. The proof is
given in the Appendix \ref{app:say}. 

\begin{lemma} \label{lm:qd2}
Assume $\psi(\theta)<\infty$ for some $\theta>0$ and $\psi_N(\theta)$ is
continuously differentiable throughout $\mathbb{R}$. We also assume the CDF
of $V_n$ is continuous and $EV_n^q<\infty$ for any $q>0$. Then
$$E[\sum_{i=1}^{N}(r_i+z)]=o(s^\delta)$$
for any $\delta>0$.
\end{lemma}

This concludes the proof of the light tailed case. We next extend the theorem
to the heavy-tailed case. We prove it by drawing connection to the truncated
system. Here we delicately choose the truncation parameter $b$ so that the
truncated system still operating the QD regime. More specifically, we choose
$b$ such that
\[
\int_{b}^{\infty} \bar G(x) dx < 1/\rho-1.
\]
This can be achieved since $EX_{n}=\int_{0}^{\infty} \bar G(x)dx<\infty$. Then
for fixed such $b$ we have
\[
\rho_{b}=\frac{E[V_{n}]}{E[X_{n}\wedge b]}=\frac{EV_{n}}{EX_{n} - \int
_{b}^{\infty}\bar G(x) dx}<1
\]
and
\[
E_{\pi}^{s} \tau(b)=o(s^{\delta})
\]
for any $\delta>0$, where $\tau(b)$ denote the coalescence time in the
truncated system.

We next prove by contradiction that the coalescence in the truncated system
implies the coalescence in the original system with the same amount of
information simulated. Recall that $\tau(b)$ is a random time satisfying that
the system has less than $s$ customers at $\tau(b)$. The maximum remaining
service time among all customers in the system at time $\tau$ is denoted as
$R(\tau(b))$. $R(\tau(b)) \leq|\tau(b)|$ and during $R(\tau(b))$ unites of
time from $\tau(b)$ on the system always has less than $s$ customers. We can
look for $\tau(b)$ at departure times of customers. We assume the process
$Q(t,y)$ is right continuous with left limit, so customers departure at time
$t$ will not counted in $Q(t,0)$. Suppose $\tau(b)$ equals to the departure
time of the $n$-th customer. Then every customer arriving between $\tau(b)$
and $\tau(b)+R(\tau(b))$ sees strictly less than $s$ customers (excluding
himself) when he enters the system. We set $\tau$ equal to the departure time
of the $n$-th customer in the original system and $R(\tau)$ by definition
equals to the maximum remaining service time among all customers in the system
at time $\tau$. We have $R(\tau) \leq R(\tau(b))$. We claim that every
customer arriving between $\tau$ and $\tau+R(\tau)$ must see less than $s$
customers (excluding himself) when he enters the system. Suppose this is not
the case. Then there exist a customer $m$, $1\leq m\leq n$ who arrives between
$\tau$ and $\tau+R(\tau)$ and finds at least $s$ customers in the system
already. The customer with the same index $m$ must have arrived between
$\tau(b)$ and $\tau(b)+R(\tau(b))$ in the truncated system and $Q(A_{m}(b)-)
\geq Q(A_{m}-) \geq s$. We get a contradiction. Therefore, we must have seen
the coalescence in the original system as well with the same amount of
information simulated.

\subsubsection{Proof of Theorem \ref{th:qed}.}

\label{sec:prof:qed}

For QED regime, by ``well-behaved", we mean that the process does not deviate $C\sqrt s$, for some $C>0$, from its fluid limit.
The following lemma states that the probability of both being well-behaved and not well-behaved are bounded away from zero. 

\begin{lemma} \label{lm:qed1}
Fix any $\eta>0$. Let $\nu(y)=(\int_{y}^{\infty} \bar F(u) du)^{1/(2+\eta)}$%
. Assume $EX_n^2<\infty$ and $EV_n^{q}<\infty$
for any $q>0$. Then for any large enough $C$, there exists $\zeta_1(C)>0$
and $\zeta_2(C)>0$, such that
\begin{equation} \label{eq:5}
P(\bar Q(t,y) \leq \lambda s \int_{y}^{t+y} \bar F(u) du + C \sqrt s \nu(y) \mbox{ for all } t\in[0,1], y \in [0, \infty)) \geq \zeta_1(C)
\end{equation}
and
\begin{equation} \label{eq:6}
P(\bar Q(t,y) > \lambda s \int_{y}^{t+y} \bar F(u) du + C \sqrt s \nu(y) \mbox{ for some } t\in[0,1], y \in [0, \infty)) \geq \zeta_2(C).
\end{equation}
\end{lemma}

The proof of Lemma \ref{lm:qed1} follows form the proof of
Lemma 9 in \cite{BlanLam:2012}. Our case is actually simpler,
as we are dealing with a one sided bound (upper bound) only as appose to the
two sided bound in \cite{BlanLam:2012}. This simplification
allows us to remove the light-tail assumption on interarrival time
distribution required in \cite{BlanLam:2012}. We shall only
briefly outline the procedure here.

For Inequality (\ref{eq:5}), the idea is to consider the diffusion limit of
$Q(t,y)$ as a two dimensional Gaussian random field \cite{WhitPang:2010}, and then invoke Borell-TIS inequality \cite{Adl:1990}. 

Assume $EX_{n}^{2} <\infty$,
$EV_{n}<\infty$ and the CDF of $V_{n}$ is continuous. Pang and Whitt
\cite{WhitPang:2010} has proved that for the $GI/GI/\infty$ queue with any
given initial age $E(0)$,
\[
\frac{\bar Q(t,y)- \lambda s\int_{t}^{t+y} \bar F(u)du}{\sqrt s} \Rightarrow
R(t,y) \mbox{ in } D_{D[0,\infty)}[0,\infty),
\]
where $R(t,y)=R_{1}(t,y)+R_{2}(t,y)$ is a Gaussian random field with
$R_{1}(t,y)=\lambda\int_{0}^{t} \int_{0}^{\infty}I(u+x>t+y)dK(u,x)$
and
$R_{2}(t,y)=\lambda c_{a}^{2}\int_{0}^{t} \bar F(t+y-u)dB(u)$,
where $K(u,x)=W(\lambda u, F(x))-F(x)W(\lambda u, 1)$ in which $W(\cdot
,\cdot)$ is a standard Brownian sheet on $[0,\infty)\times[0,1]$ and
$B(\cdot)$ is a standard Brownian motion independent of $W(\cdot,\cdot)$. The
constant $c_{a} $ is coefficient of variation of the interarrival times, i.e.
$c_{a}=\sqrt{\mbox{Var}(X_{n})}/EX_{n}$. We denote
\[
\tilde R_{i}(t,y)=\frac{R_{i}(t,y)}{v(y)}%
\]
and define the d-metric (a pseudo-metric)
\[
d_{i}((t,y),(t^{\prime},y^{\prime}))=E[(\tilde R_{1}(t,y)-\tilde
R_{2}(t^{\prime},y^{\prime}))^{2}]
\]
for $i=1,2$.\newline We then invoke the Borell-TIS inequality. We shall skip
the verification of the conditions for such invocation here as it is tedious
and detailedly proved in \cite{BlanLam:2012}. Let
$S=[0,1]\times[0,\infty)$. it is shown in \cite{BlanLam:2012} that, 
there exist constants $M_{i,1}>0$ and $M_{i,2}>0$, such
that $E[\sup_{S} \tilde R_{i}(t,y)] \leq M_{i,1} <\infty$ and $\sup_{S}
E[\tilde R_{i}(t,y)^{2}]\leq M_{i,2} < \infty$. And for $C_{i} \geq E[\sup_{S}
\tilde R_{i}(t,y)]$,
\[
P(\sup_{S} \tilde R_{i}(t,y) \geq C_{i}) \leq\exp\{-\frac{1}{2\sup_{S}E[\tilde
R_{i}(t,y)^{2}]}(C_{i}-E[\sup_{S} \tilde R_{i}(t,y)])^{2}\}
\]
for $i=1,2$. \newline Let $C \geq2\max\{E[\sup_{S} \tilde R_{1}(t,y)],
E[\sup_{S} \tilde R_{2}(t,y)]\}$. Then
\begin{align*}
&  P(R(t,y) \leq C\nu(y) \mbox{ for all } t \in[0,1], y\in[0,\infty))\\
&  \geq P(\sup_{S} \tilde R_{1}(t,y) + \sup_{S} \tilde R_{2}(t,y) \leq C)\\
&  \geq P(\sup_{S} \tilde R_{1}(t,y) \leq\frac{C}{2})P(\sup_{S}\tilde
R_{2}(t,y) \leq\frac{C}{2}) >0.
\end{align*}

Let $X_{0}$ denote the interarrival time of the first customer and $V_{0}$
denote its service time. We also denote $\bar Q^{0}(t,y)$ as an independent
infinite server process starting empty and with $E(0)=0$. Then for $s$ large
enough, we have
\begin{align*}
&  P(\bar Q(t,y) \leq\lambda s \int_{y}^{t+y} \bar F(u)du + C \sqrt s \nu(y)
\mbox{ for all } t\in[0,1], y\in[0,\infty))\\
&  =P(\bar Q^{0}(t-X_{0},y) + 1\{V_{0} > t+y\} \leq\lambda s \int_{y}^{t+y}
\bar F(u)du + C \sqrt s \nu(y)\\
& \mbox{ for all } t\in[X_{0},1], y \in [0,\infty))\\
&  \geq P(\bar Q^{0}(t,y)+1\{V_{0}>t+X_{0}+y\} \leq\lambda s \int_{y}%
^{t+X_{0}+y} \bar F(u)du + C\sqrt s\nu(y)\\
& \mbox{ for all } t\in[0,1-X_{0}],
y\in[0,\infty))\\
&  \geq P(\bar Q^{0}(t,y) \leq\lambda s \int_{y}^{t+y} \bar F(u)du + C \sqrt s
\nu(y) \mbox{ for all } t\in[0,1], y\in[0,\infty))\\
&  = P(\frac{\bar Q^{0}(t,y)-\lambda s\int_{y}^{t+y} \bar F(u) du}{\sqrt s}
\leq C \nu(y) \mbox{ for all } t\in[0,1], y\in[0,\infty)).
\end{align*}
It is easy to check that the set $\{f: |f(t,y)| \leq C\nu(y) \mbox{ for all }
t\in[0,1], y\in[0,\infty)\}$ is a continuity set, thus by the Functional
Central Limit Theorem result in \cite{WhitPang:2010},
we have
\begin{align*}
&  P(\frac{\bar Q^{0}(t,y)-\lambda s\int_{y}^{t+y} \bar F(u) du}{\sqrt s} \leq
C \nu(y) \mbox{ for all } t\in[0,1], y\in[0,\infty))\\
&  \rightarrow P(R(t,y) \leq C\nu(y) \mbox{ for all } t\in[0,1], y\in
[0,\infty))>0.
\end{align*}

Inequality (\ref{eq:6}) is easy to prove as we can always isolate a point
$(t^{*},y^{*})$ inside $S$. The projection of the process on that point posses
Gaussian distribution. More specifically,
\begin{align*}
&  P(\bar Q(t,y) > \lambda s \int_{y}^{t+y} \bar F(u) du + C \sqrt s \nu(y)
\mbox{ for some } t\in[0,1], y \in[0, \infty))\\
&  \geq P(\bar Q(t^{*}, y^{*}) > \lambda s \int_{y^{*}}^{t^{*}+y^{*}}\bar
F(u)du+C\sqrt s\nu(y^{*}))\\
&  = P(\frac{\bar Q(t^{*}, y^{*})-\lambda s \int_{y^{*}}^{t^{*}+y^{*}}\bar
F(u)du}{\sqrt s} > C\nu(y^{*})),
\end{align*}
and by Fatou's lemma
\[
\liminf_{s\rightarrow\infty} P(\frac{\bar Q(t^{*}, y^{*})-\lambda s
\int_{y^{*}}^{t^{*}+y^{*}}\bar F(u)du}{\sqrt s} > C\nu(y^{*})) \geq
P(R(t^{*},y^{*})> C\nu(y^{*}))>0.
\]

Let $m(s)=\inf\{y: C\sqrt s(v(y)+\int_{y}^{\infty}v(s)ds) \leq\frac{1}{2}\}$.
Following the same construction as for the QD regime, we
define the sequence of random times $\Xi_{i}$'s as follows: $\Xi_{0}:=0$.
Given $\Xi_{i-1}$ for $i=1,2,\cdots$,
\begin{align*}
r_{i}  &  =\inf\{k: k \geq R(\Xi_{i-1}), k=1,2,...),\\
z  &  = \inf\{k: k \geq m(s), k=1,2,...\},\\
\Xi_{i}  &  = \Xi_{i-1}+r_{i}+z.
\end{align*}
We introduce a Bernoulli random variable $\xi_{i}$ with $\xi_{i}=1$ if and
only if
\[
\bar Q_{\Xi_{i-1}+(k-1)t_{0}}(t,y) \leq\lambda s\int_{y}^{t+y} \bar F(u)du + C
\sqrt s\nu(y)
\]
for all $t \in[0, 1]$, $y\in[0,\infty)$ and every $k=1,2,...,r_{i}+z$.

We next show that $\xi_{i}=1$ implies that
$\tau^{*}$ is reached before $\Xi_{i}$. Since $r_{i} \geq R(\Xi_{i-1})$, all
the customers at time $\Xi_{i-1}+r_{i}$ will be those arrive after $\Xi_{i}$.
Thus we have $\xi_{i}=1$ implies that
\begin{align*}
Q(\Xi_{i-1}+r_{i}, y)  &  \leq\sum_{k=1}^{r_{i}} \{ \lambda s\int
_{(k-1)t_{0}+y}^{kt_{0}+y}\bar F(u)du + C \sqrt s \nu((k-1)+y))\}\\
&  \leq\lambda s \int_{y}^{\infty}\bar F(u)du + C \sqrt s (\nu(y)+\int
_{y}^{\infty}\nu(u)du).
\end{align*}
As $\int_{y}^{\infty} \bar F(u)du$ decays faster than $\nu(y)$ as $y$ grows
large, for $s$ large enough, we have
\[
R(\Xi_{i-1}+r_{i})<m(s).
\]
Likewise for every $t\in(k-1,k]$ and $k=1,2,\cdots,z$,
\[
Q(\Xi_{i-1}+r_{i}+t, y) \leq\lambda s \int_{y}^{\infty}\bar F(u)du + C \sqrt s
( \nu(y)+\int_{y}^{\infty}\nu(u)du).
\]
Thus when $\beta>C(\nu(0)+\int_{0}^{\infty} \nu(u)du)$, we have
\[
Q(\Xi_{i-1}+r_{i}+t, 0) \leq s+C(\nu(0) + \int_{0}^{\infty} \nu(u)du) \sqrt s
\leq s+\beta \sqrt s
\]
for $t \in[0, R(\Xi_{i-1}+r_{i})]$.\newline Now let $N=\inf\{i\geq1: \xi
_{i}=1\}$. Then
\[
E\tau^{*} \leq E[\sum_{i=1}^{N}(r_{i}+z)].
\]

We now show a bound for $E\sum_{i=1}^{N}(r_{i}+z)$. The proof is given in the Appendix \ref{app:say}. 

\begin{lemma} \label{lm:qed2}
Assume $EX_n^2<\infty$ and $EV_n^q<\infty$ for any $q>0$. Then
$$\log E[\sum_{i=1}^{N}(r_i+z)]=o(s^\delta)$$
for any $\delta>0$.
\end{lemma}

Notice that our proof of Theorem \ref{th:qed} only requires the existence of
the second moment of the interarrival time distribution. We thus conclude the proof of Theorem \ref{th:qed}.

%If your paper includes appendices, then precede the first of them by the command
\appendix
%and then carry on using the \section and \subsection commands, as above.

\section{Proof of Proposition \ref{prop}} \label{app:prop}
By Chebyshev's inequality,
$$P(A_{n+1} < n(\mu-\epsilon)) \leq E[\exp(\theta(n(\mu-\epsilon)-A_{n+1}))] \leq \exp(-n(-\theta(\mu-\epsilon)-\psi(-\theta)))$$
for any $\theta \geq 0$.\\
Let
$$I(-\epsilon)=\max_{\theta \geq 0}\{-\theta(\mu-\epsilon)-\psi(-\theta)\}.$$
As $\psi(0)=0$, $\psi'(0)=\mu$ and $\psi''(0)=Var(X)>0$, $I(-\epsilon)>0$. Then
$$P(A_{n+1}<n(\mu-\epsilon))\leq\exp(-nI(-\epsilon))$$
and
$$\sum_{n=1}^{\infty} P(A_{n+1} < n(\mu-\epsilon)) \leq \frac{\exp(-I(-\epsilon))}{1-\exp(-I(-\epsilon))} < \infty.$$
By Borel-Cantelli lemma, $\{A_{n+1} \geq n(\mu-\epsilon)\}$ eventually almost surely.\\
Similarly and independently we have
\begin{eqnarray*}
\sum_{n=1}^{\infty} P(\left\vert V_{n+1} \right\vert>(n(\mu-\epsilon))^{\alpha}) &=& \sum_{n=1}^{\infty} P(\left\vert V_1 \right\vert^{1/\alpha}>n(\mu-\epsilon))\\
&\leq& \frac{1}{\mu-\epsilon}\int_{0}^{\infty} P(\left\vert V_1 \right\vert^{1/\alpha}>\nu)d\nu < \infty.
\end{eqnarray*}
Thus, again by Borel-Cantelli lemma, $\{\left\vert V_{n+1} \right\vert \leq (n(\mu-\epsilon))^{\alpha}\}$ eventually almost surely. Therefore,
$P(\kappa<\infty)=1$.

\section{Proof of Lemma \ref{lm:qd2} and Lemma \ref{lm:qed2}} \label{app:say}

\begin{lemma} \label{lm:max}
If $EV_n^q<\infty$ for any $q>0$, then for any fixed $p>0$,
\begin{equation*}
E[(\max_{k=1,2,...n} V_k)^p]=o(n^\delta)
\end{equation*}
for any $\delta>0$.
\end{lemma}

\begin{proof}
For any fixed $\delta>0$ we can find $\delta' \in (0,\delta)$. Let $q=1/\delta'+p$. By Chebyshev's inequality we have
$$\bar F(u) \leq \frac{EV^{q}}{u^q}.$$
Let $\bar F_n(u)=P(\max_{k=1,2,...,n} V_k > u)$ then
\begin{eqnarray*}
E[(\max_{k=1,2,...n} V_k)^p] &=& p\int_{0}^{\infty} u^{p-1} \bar F_n(u) du\\
&\leq& n^{1/(q-p)}+np\int_{n^{1/(q-p)}}^{\infty} u^{p-1}\bar F(u) du\\
&\leq& n^{1/(q-p)} + np\int_{n^{1/(q-p)}}^{\infty} \frac{EV^{q}}{u^{q-p+1} }du\\
&=& n^{\delta'} + \frac{p}{q-p} EV^{q}.
\end{eqnarray*}
\end{proof}

\begin{eqnarray*}
E[\sum_{i=1}^{N}(r_i+z)] &=& E[\sum_{i=1}^{\infty} (r_i+z)I\{N \geq i\}]\\
&\leq& \sum_{i=1}^{\infty}E[(r_i+z)^2]^{1/2}P(N \geq i)^{1/2} \mbox{ by Holder's inequality}.
\end{eqnarray*}

\begin{lemma}
If $EX_n<\infty$ and $EV_n^{q}<\infty$ for any $q>0$, then for any $p\geq1$ we have
$$E[(r_i+z)^p]^{1/p}=o(s^\delta)$$
for any $\delta>0$.
\end{lemma}

\begin{proof}
By Minkowski inequality
$$E[(r_i+z)^p]^{1/p}\leq E[r_i^p]^{1/p}+z.$$
Using similar argument as in the proof of Lemma \ref{lm:max}, we can show that $l(s)=o(s^{\delta})$ for any $\delta>0$, thus $z=o(s^{\delta})$ for any $\delta>0$.\\
For fixed $\delta>0$, we can find $\delta' \in (0, p\delta/(1+p\delta))$, such that
\begin{eqnarray*}
E[r_i^p] &\leq& E[E[ (\max_{k=1,...,N_s(\Xi_{i-1})-N_s(\Xi_{i-2})} V_k )^p|N_s(\Xi_{i-1})-N_s(\Xi_{i-2})]]\\
&\leq& CE[(N_s(\Xi_{i-1})-N_s(\Xi_{i-2}))^{\delta'}]   \mbox{  Lemma \ref{lm:max}}\\
&\leq& C(E[N_s(\Xi_{i-1})-N_s(\Xi_{i-2})])^{\delta'}    \mbox{  Jensen's inequality for concave function} \\
&\leq& C\tilde\lambda^{\delta'} s^{\delta'} E[r_{i-1}+z]^{\delta'}  \mbox { Key Renewal Theorem}.
\end{eqnarray*}
Let $w_i=r_i+z$ for $i=1,2,\cdots$. As $z$ it is a constant that only depends on $s$ and $z=o(s^{\delta'})$, then
$$Ew_i \geq z \geq 1$$
and
$$Ew_i = Er_i+z \leq C \tilde \lambda^{\delta'} s^{\delta'} (Ew_{i-1})^{\delta'} + z \leq \tilde C s^{\delta'} (Ew_{i-1})^{\delta'}$$
where $\tilde C=C\tilde\lambda^{\delta'}+1$.\\
As $E[r_1^p]=E_\pi[R(0)^p]=o(s^{\delta'})$. By iteration we have
$$Ew_i \leq \tilde C^{1/(1-\delta^{\prime})} s^{\delta'/(1-\delta^{\prime})}$$
for $i=1,2,\cdots$.\\
Thus $Er_i^p=o(s^{p\delta})$ and $E[(r_i+z)^p]^{1/p}=o(s^\delta)$.
\end{proof}

\begin{proof}[Proof of Lemma \ref{lm:qd2}]
We first notice that $P(\xi_i=0)\leq E[w_1] \exp(-sI^*(\delta)+o(s))$ by Lemma \ref{lm:dp1}.
$$P(N \geq 1)=1.$$
$$P(N \geq 2)=P(\xi_1=0) \leq E[w_1] \exp(-sI^*(\delta)+o(s))$$
Recall that $w_i=r_i+z$ for $i=1,2,\cdots$.
\begin{eqnarray*}
P(N\geq3) &=& P(N\geq 1)P(N\geq 3|N\geq 2)\\
&=& P(\xi_1=0) P(\xi_2=0|\xi_1=0)\\
&\leq& P(\xi_1=0)E[w_2|\xi_1=0] \exp(-sI^*(\delta)+o(s))\\
&\leq& E[w_1]E[w_2|\xi_1=0] \exp(-2sI^*(\delta)+o(s)).
\end{eqnarray*}
We next prove that $E[w_2|\xi_1=0]=\exp(o(s))$. Notice that $P(\xi_i=0) \geq \exp(-sI^*(\delta)+o(s))$ by Lemma \ref{lm:dp1}. Then for any $p>0$, $q>0$ and $1/p+1/q=1$,
\begin{eqnarray*}
E[w_2|\xi_1=0] &=& \frac{E[w_2I\{\xi_1=0\}]}{P(\xi_1=0)}\\
&\leq& \frac{E[w_2^p]^{1/p} P(\xi_1=0)^{1/q}}{P(\xi_1=0)} \mbox{ Holder's inequality }\\
&\leq& E[w_2^p]^{1/p}E[w_1]^{1/q}\exp(\frac{1}{p}sI^*(\delta)+o(s)),
\end{eqnarray*}
thus
$$\frac{1}{s} \log E[w_2|\xi_1=0] \leq \frac{1}{s}(\frac{1}{p}\log E[w_2^p] +\frac{1}{q} \log E[w_1]+o(s)) + \frac{1}{p} I^*(\delta).$$
By sending $p$ to infinity, we have $E[w_2|\xi_1=0] = \exp(o(s))$.\\
Similarly by iteration,
$$P(N\geq k)=\exp(-ksI^*(\delta)+o(s))$$
for $k = 4,5,\cdots$.\\
Then
$\sum_{i=1}^{\infty} P(N \geq i)^{1/2} =O(1).$
As $E[\sum_{i=1}^{N}(r_i+z)] \leq \sum_{i=1}^{\infty} E[(r_i+z)^2]^{1/2} P(N \geq i)^{1/2}$ and $E[(r_i+z)^2]^{1/2}=o(s^\delta)$ for any $\delta>0$, we have
$$E[\sum_{i=1}^{N}(r_i+z)]=o(s^\delta)$$
for any $\delta>0$.
\end{proof}

\begin{proof}[Proof of Lemma \ref{lm:qed2}]
$$P(N\geq1)=1.$$
\begin{eqnarray*}
P(N\geq2) &=& P(\xi_1=0)\\
&\leq& 1- E[\zeta_1(C)^{w_1}] \mbox{ Lemma \ref{lm:qed1}}\\
&\leq& 1- \zeta_1(C)^{E[w_1]}  \mbox{ Jensen's inequality}\\
&=& 1-b\exp(-o(s^\delta)).
\end{eqnarray*}
Moreover
\begin{eqnarray*}
P(N \geq 3) &=& P(N>2|N>1) P(N>1)\\
&=& P(\xi_2=0|\xi_1=0) P(\xi_1=0)\\
&\leq& E[1-\zeta_1(C)^{w_2}|\xi_1=0] P(\xi_1=0)\\
&\leq& (1-\zeta_1(C)^{E[w_2|\xi_1=0]})P(\xi_1=0).
\end{eqnarray*}
We next show that $E[w_2|\xi_1=0]=o(s^\delta)$ for any $\delta>0$. Notice that $P(\xi_i=0) \geq \zeta_2(C)$ by Lemma \ref{lm:qed1}, then
$$E[w_2|\xi_1=0] =\frac{E[w_2 I\{\xi_1=0\}]}{P(\xi_1=0)} \leq \frac{Ew_2}{\zeta_2(C)}$$
Similarly by iteration we have
$$P(N \geq k) \leq (1-b\exp(-o(s^\delta)))^k$$
for any $\delta>0$ and $k=4,5,\cdots$.\\
Then
$$\log \sum_{i=1}^{\infty} P(N \geq i)^{1/2} =o(s^{\delta})$$
for any $\delta>0$.\\
As $E[\sum_{i=1}^{N}(r_i+z)] \leq \sum_{i=1}^{\infty} E[(r_i+z)^2]^{1/2} P(N \geq i)^{1/2}$ and $E[(r_i+z)^2]^{1/2}=o(s^\delta)$ for any $\delta>0$, we have
$$\log E[\sum_{i=1}^{N}(r_i+z)]=o(s^{\delta})$$
for any $\delta>0$.\\
\end{proof}

\noindent{\bf\large Acknowledgement:}
NSF support from grants CMMI-0846816 and 1069064 is gratefully acknowledged.
% Place the text of your acknowledgements after the \acks command.
% \acks generates the heading "Acknowledgements".
% If you wish to make only one acknowledgement, use \ack.
% \ack generates the heading "Acknowledgement".

% Reference list
%
% References should be in the following form (or the BibTeX file
% apt.bst should be used):
%
% For a journal:
% Surname, Initial (year). Title of paper. {\em Journal title}
% {\bf Vol,} page--range.
%
% For a book:
% Surname, Initial (year). {\em Book title}. Publisher, Address.
%
% Note the following example of a reference list.

\bibliographystyle{apt}
\bibliography{exact}

\end{document}